\documentclass{amsproc}
\usepackage{graphicx}
\usepackage{amssymb}
\usepackage{epstopdf}
\usepackage{amsmath,amsthm,amscd,amssymb}
\usepackage{latexsym}
\usepackage[colorlinks,citecolor=red,pagebackref,hypertexnames=false]{hyperref}
\usepackage{geometry}                
\numberwithin{equation}{section}

\theoremstyle{plain}
\newtheorem{theorem}{Theorem}[section]
\newtheorem{lemma}[theorem]{Lemma}
\newtheorem{corollary}[theorem]{Corollary}

\newtheorem{conjecture}[theorem]{Conjecture}

\theoremstyle{definition}
\newtheorem{definition}[theorem]{Definition}

\newtheorem{case[theorem]}{Case}

\newcommand{\atxt}{}

\def\R{\mathbb R}
\def\sd{\mathbb S^{d-1}}
\def\implies{\Longrightarrow}
\def\O{\mathcal O}

\def\zero{\mathbf 0}
\def\be{\begin{equation}}
\def\ee{\end{equation}}
\def\bes{\begin{eqnarray*}}
\def\ees{\end{eqnarray*}}
\def\st{\sqrt{3}}
\def\I{\mathbb I}

\theoremstyle{remark}
\newtheorem{remark}[theorem]{Remark}

\numberwithin{equation}{section}

\begin{document}

\title[Multilinear generalized Radon transforms and point configurations] {\parbox{14cm}{\centering{Multilinear generalized Radon transforms  \\ and point configurations}}}


\author{Loukas Grafakos, Allan Greenleaf, Alex Iosevich and Eyvindur Palsson}

\address{Department of Mathematics \\ University of Missouri\\ Columbia, MO 65211}
\email{grafakosl@missouri.edu}

\address{Department of Mathematics \\ University of Rochester\\ Rochester, NY 14627}
\email{allan@math.rochester.edu}

\address{Department of Mathematics \\ University of Rochester\\ Rochester, NY 14627}
\email{iosevich@math.rochester.edu}

\address{Department of Mathematics \\ University of Rochester\\ Rochester, NY 14627}
\email{palsson@math.rochester.edu}
\thanks{The authors were partially supported by NSF grants   DMS 0900946, DMS-0853892 and DMS-1045404.}

\begin{abstract}

We study multilinear generalized Radon transforms using a graph-theoretic paradigm that includes the widely studied linear case. 
These provide a general mechanism to study Falconer-type problems involving $(k+1)$-point configurations in geometric measure theory,  with $k \ge 2$, including the  distribution of simplices, volumes and angles determined by the points of
fractal subsets $E \subset {\Bbb R}^d$, $d \ge 2$. 
If $T_k(E)$ denotes the set of noncongruent $(k+1)$-point configurations determined by  $E$, we show that  if the Hausdorff dimension of $E$ is greater than $d-\frac{d-1}{2k}$, then the ${k+1 \choose 2}$-dimensional Lebesgue measure of $T_k(E)$ is positive. This complements  previous work on the Falconer conjecture (\cite{Erd05} and the references there), as well as work on finite point configurations  \cite{EHI11,GI10}. We also give applications to 
Erd\"os-type problems in discrete geometry and a fractal regular value theorem, providing a multilinear framework for the results in \cite{EIT11}. 
\end{abstract}

\maketitle


\section{Introduction}

\vskip.125in

{\atxt Linear generalized Radon transforms  are operators of the form
\begin{equation} \label{phongsteindef} 
{\mathcal R}f(x)=\int_{\{\Phi(x,y)=\vec{t}\}} f(y) \Psi(x,y) d\sigma_{\vec{t}}^x(y), 
\end{equation}
where $\Phi:{\Bbb R}^d \times {\Bbb R}^d \to {\Bbb R}^m$, $d \ge 2$, is a family of smooth defining  functions, $\Psi$ is a smooth cut-off function,  and $d\sigma_{\vec{t}}^x(y)$ on $\{y: \Phi(x,y)=\vec{t}\}$ is induced from the Leray measure on the incidence relation  $\Sigma_{\vec{t}}:=\{(x,y): \Phi(x,y)=\vec{t}\}$. All of these objects, and the operator $\mathcal R$, vary smoothly as $\vec{t}$ varies over regular values of $\Phi$.
Operators  on this level of generality were introduced by Guillemin and Sternberg \cite{GuSt} and Phong and Stein \cite{PS83}, building on earlier work of Gelfand and his collaborators, and Helgason. Beyond their original role in integral geometry, generalized Radon transforms and their singular variants have since become ubiquitous in harmonic analysis, partial differential equations and related areas \cite{PS86,PS86II,PS89,PS91,GS02}, and more recently as  tools to study geometric and combinatorial problems 
\cite{Fal86, Erd05,IJL09}. }

Model cases of generalized Radon transforms have been present in the literature for a long time. 
An example of particular importance to this paper is the spherical averaging operator, 
\begin{equation} \label{spherical} A^d_1f(x)=\int_{S^{d-1}} f(x-y) d\sigma(y), \end{equation} 
$d \ge 2$, where $d\sigma$ is the Lebesgue measure on the unit sphere $S^{d-1}$. 
It was proved in \cite{Litt73,Str70} that 
\begin{equation} \label{sphericalsobolev} A^d_1: L^2({\Bbb R}^d) \rightarrow L^2_{\frac{d-1}{2}}({\Bbb R}^d),\end{equation} 
where $L^2_s({\Bbb R}^d)$ is the standard Sobolev space consisting of $L^2({\Bbb R}^d)$ functions with $s$ generalized derivatives in $L^2({\Bbb R}^d)$, and that 
\begin{equation} \label{sphericalpq} A^d_1: L^p({\Bbb R}^d) \rightarrow L^q({\Bbb R}^d) \end{equation} if and only if $\big(\frac{1}{p}, \frac{1}{q} \big)$ belongs to the closed triangle with the endpoints $(0,0)$, $(1,1)$ and $\big(\frac{d}{d+1}, \frac{1}{d+1} \big).$ This led to  ongoing studies of $L^p$-improving properties of measures; see, e.g., \cite{C98,O08,TW03} and the references  there. Among many other applications, $A^d_1$ is closely related to the fundamental solution of  the wave equation; see \cite{J,So93} and the references  there. 

The purpose of this paper is to introduce a class of multilinear generalized Radon transforms and  apply them to several problems in harmonic analysis,  geometric measure theory and discrete geometry. Before initiating this study, let us reinterpret the definition  of the (linear) generalized Radon transforms  in terms of a straight forward graph-theoretic paradigm.
Given  $\Phi:{\Bbb R}^d\times {\Bbb R}^d \to {\Bbb R}^m$ as in (\ref{phongsteindef}), define a (directed)
graph whose vertices are points in  $ {\Bbb R}^d$, by saying  that two vertices,  $x, y \in \Bbb R^d$, are connected by an \emph{edge}, or $x \sim y$,  if $\Phi(x,y)=\vec{t}$. Recall that the adjacency operator for a finite graph $G$ is defined by 
$$ Af(x)=\sum_{x \sim y} f(y).$$ Thus, the generalized Radon transform 
(\ref{phongsteindef}) may be viewed as a continuous analogue of the adjacency operator for the infinite directed graph defined by the pair $(\Phi,\vec{t})$. 

The graph-theoretic perspective on generalized Radon transforms is implicit in the work of Falconer \cite{Fal86} and subsequent efforts on the Falconer distance problem. The question there is to determine how large the Hausdorff dimension of a subset $E\subset{\Bbb R}^d$, $d \ge 2$, needs to be to ensure that the Lebesgue measure of the set of pairwise distances is positive. This means that one must show that a given distance cannot arise too often, either in a point-wise or average sense. 
If we view the points of the ambient set as vertices of a graph and connect two vertices by an edge if they are separated by a given fixed distance, then the problem is to obtain a suitable bound on the distribution of edges of this graph. This naturally leads one to the examination of the Sobolev bounds for the adjacency operator which, in this case, turns out to be the spherical averaging operator defined in (\ref{spherical}). 

The distance problem can be viewed as a geometric problem on two-point configurations  in subsets of ${\Bbb R}^d$; in this paper, we introduce  $k$-multilinear variants of  generalized Radon transforms which can be used to study $(k+1)$-point configurations, $k\ge 2$. 

\subsection{Definition of a multilinear generalized Radon transform}\label{multidefsubsec} The graph-theoretic point of view leads  naturally to the 
definition of multilinear generalized Radon transforms. Consider a (directed) hyper-graph whose vertices are points in  $ {\Bbb R}^d$. For  $1 \le k \le d$, 
let  $\Phi:(\mathbb R^d)^{k+1}\rightarrow\mathbb R^m$ and $\vec{t}\in\mathbb R^m$. Then,
we say that  the ordered $(k+1)$-tuple of vectors $(x^1, x^2, \dots, x^{k+1})$ are connected by a \emph{hyper-edge} if $\Phi(x^1,\dots,x^{k+1})=\vec t$. Letting 
$$\Sigma_{\vec{t}} =\big\{(x^1, \dots , x^{k+1}) \in (\mathbb R^d)^{k+1}: \Phi(x^1,\dots,x^{k+1})=
 \vec t\, \big\},$$ 
 the continuous variant of the adjacency operator for this hyper-graph is given by


\begin{equation} \label{multidef} {\mathcal R}^d_k(f_1, \dots, f_k)(x^{k+1})=
\int_{\Sigma_{\vec{t}}^{x^{k+1}}}\quad
 \prod_{j=1}^k f_j(x^j)\, d\sigma_{\vec{t}}^{x^{k+1}}(x^1, \dots, x^k), \end{equation} 
where $d\sigma_{\vec{t}}^{x^{k+1}}(x^1, \dots, x^{k})$ is the Leray measure on the set 
$$\Sigma_{\vec{t}}^{x^{k+1}} =\big\{(x^1, \dots, x^k) \in (\mathbb R^d)^{k}: \Phi(x^1,\dots, x^{k+1})=
 \vec t\, \big\}.$$ 
This may be modified  in an inessential way by multiplying $d\sigma_{\vec{t}}^{x^{k+1}}$ by a smooth cut-off function $\psi$. 

A model $k$-linear generalized Radon transform is the multilinear  analogue of the linear spherical averaging operator.
Taking advantage of its translation invariance, we define this as
\begin{equation} A^d_k(f_1, \dots, f_k)(x)=\int \dots \int \prod_{j=1}^k f_j(x-u^j) dM_k^d(u^1, \dots, u^k), \end{equation} 
where,  for $k \leq d$,  $dM_k^d$ is the Leray measure on the set 
$$ \Sigma_d^k=\{(u^1, \dots, u^k) \in S^{d-1} \times \dots \times S^{d-1}: |u^i-u^j|=1; 1 \leq i<j \leq k \}.$$ 
Just as the spherical averaging operator arose naturally in Falconer's and subsequent investigations of the distance problem, the $k$-linear operator $A_k^d$ arises naturally when considering the $(k+1)$-point configuration problem investigated in \cite{GI10,EHI11,IMP11}. The question there was to determine how large the Hausdorff dimension of a subset of ${\Bbb R}^d$ needs to be to ensure that the ${k+1 \choose 2}$-dimensional Lebesgue measure of non-congruent $k$-simplices (i.e., $(k+1)$-point configurations) is positive. 

The operator $A_k^d$ and translation-invariant multilinear generalized Radon transforms in general can be put into context  via a result of the first listed author and Soria \cite{GS10} (stated there for bilinear convolution  operators, but easily extended to multilinear ones), and an immediate consequence. 
\begin{theorem} \label{gs} (Grafakos and Soria 2010) Let $\mu$ be a non-negative Borel measure on $(\Bbb R^d)^k$ and set
$$
T_\mu(f_1,\dots , f_k)(x) =\int_{\mathbb R^d} \cdots \int_{\mathbb R^d}
f_1(x-u^1)\cdots f_k(x-u^k)d\mu(u^1,\dots, u^k).
$$
Suppose that $\frac{1}{p_{1}}+\dots +\frac{1}{p_k}=\frac{1}{r} \leq 1$. Then 
$$T_{\mu}:L^{p_1}({\Bbb R}^d) \times\dots\times L^{p_k}({\Bbb R}^d) \to L^r({\Bbb R}^d)$$ if and only if $\mu$ is a finite measure. 
Furthermore, suppose that 
$\gamma_j\in\Bbb R$ are such that the  distribution 
$(I-\Delta_{u_1})^{\frac{\gamma_1}2}\cdots 
(I-\Delta_{u_k})^{\frac{\gamma_k}2}\mu$  is a Borel measure.
Then for all $p_j$ with $\frac{1}{p_1}+\cdots +\frac{1}{p_k}=\frac{1}{\tilde{r} }\le 1$,  
$$T_\mu:
 L^{p_1}_{-\gamma_1}(\mathbb R^d)\times \cdots \times 
 L^{p_k}_{-\gamma_k}(\mathbb R^d)
\to L^{\tilde{r}}(\mathbb R^d),$$
where $ L^{p}_\gamma$ denotes the inhomogeneous Sobolev space of  distributions $g$ with $(I-\Delta)^\gamma g \in L^p(\mathbb R^d)$. 
\end{theorem} 

The second part of Theorem \ref{gs}   follows easily by expressing $T_\mu$ in Fourier multiplier form 
$$
T_\mu(f_1,\dots , f_k)(x) =\int_{\mathbb R^d} \cdots \int_{\mathbb R^d}
\widehat{f_1}(\xi_1)\cdots 
\widehat{f_k}(\xi_k) \widehat{\mu}(\xi_1, \dots, \xi_k)
e^{2\pi i x\cdot(\xi_1+\cdots + \xi_k)} d\xi_1\cdots d\xi_k\, , 
$$
multiplying and dividing by $(1+4\pi^2|\xi_1|^2)^{\gamma_1/2}\cdots 
(1+4\pi^2|\xi_k|^2)^{\gamma_k/2}$, taking inverse Fourier transforms, and using the first part of the theorem.

As a special case, consider the bilinear operator in the plane, i.e. $k=d=2$, for which
\begin{equation}\label{att}
A_2^2(f,g)(x)=\sum_{\pm}\int_0^{2\pi} f(x-(\cos(\theta), \sin(\theta)) g(x-(\cos(\theta\pm\pi/3)), \sin(\theta\pm\pi/3)) 
d\theta.
\end{equation}

In the notation of Theorem \ref{gs},  the measure $\mu$ is a multiple of arc-length  on the curve
$$ \big\{(u,v) \in \mathbb R^2 \times \mathbb R^2: |u|=|v|=|u-v|=1 \big\}\subset\Bbb R^4.$$ 
One can calculate \cite{GI10} that 
$ \widehat{\mu}(\xi, \eta)= \sum_{\pm}\widehat{\sigma}(U_\pm(\xi, \eta))$, where 
$ U_\pm: {\Bbb R}^4 \to {\Bbb R}^2$  are the linear maps
$$ U_\pm(\xi, \eta)=\left( \xi_1+\frac{\eta_1}{2}\pm\eta_2 \frac{\sqrt{3}}{2}, 
\xi_2\mp\eta_1 \frac{\sqrt{3}}{2}+\frac{1}{2} \eta_2 \right).$$ 
It follows that $\widehat{\mu} \in L^{\infty}({\Bbb R}^2 \times {\Bbb R}^2)$, with no better uniform decay; however, the operator $A_2^2$ satisfies much better bounds than those implied by Theorem \ref{gs}.
In fact, we shall see in the sequel that 
\begin{equation} \label{bicircle} A_2^2: L^2_{-\frac{1}{2}}({\Bbb R}^2) \times L^2({\Bbb R}^2) 
\rightarrow L^1({\Bbb R}^2)\end{equation} and, more generally, 
\begin{equation} \label{bisphere} A_2^d: L^2_{-\frac{d-1}{2}}({\Bbb R}^d) \times L^2({\Bbb R}^d) \rightarrow L^1({\Bbb R}^d),\end{equation} 
with corresponding non-trivial bounds for the operator $A_k^d$, $k \leq d$.

In the first result of the current work,  we establish certain bounds for multilinear generalized Radon transforms. 
Later,   Theorem \ref{generalmultilinear} and its corollaries will give results for nontranslation invariant multilinear generalized Radon transforms.
However, all of our applications to continuous and discrete geometry are in fact made using the more restrictive class of translation invariant multilinear generalized Radon transforms and in this setting one can obtain  stronger results.

\begin{theorem} \label{translationinvariant} Let $T_{\mu}$ be the multilinear convolution operator 
$$ T_{\mu}(f_1, \dots, f_k)(x)=\int \dots \int f_1(x-u^1) \dots f_k(x-u^k) d\mu(u^1, \dots, u^k)$$ where $\mu$ is a nonnegative Borel measure. Suppose that 
\begin{equation}\label{MainDecay}
|\widehat{\mu}(-\xi,\xi, 0, \dots, 0)| \lesssim {(1+|\xi|)}^{-\gamma}
\end{equation}
for some $\gamma>0$. Then, for all $\gamma_1,\gamma_2 > 0$ such that $\gamma = \gamma_1 + \gamma_2$, and  acting on nonnegative functions,
$$ T_{\mu}:L^2_{-{\gamma_1}}({\Bbb R}^d) \times L^2_{-{\gamma_2}}({\Bbb R}^d) \times L^{\infty}({\Bbb R}^d) \times \dots \times L^{\infty}({\Bbb R}^d)  \to L^1(\Bbb R^d).$$ 
\end{theorem} 

Here, and throughout, the notation $X \lesssim Y$ means that there exists a constant $C>0$, independent of the variables of interest (depending on the setting), such that $X \leq CY$.

In  Theorem \ref{translationinvariant} (and Theorem \ref{generalmultilinear}  below) there is nothing special about the first two coordinates in the assumption  of Fourier decay of the measure; one could of course state both theorems more generally for any two distinct coordinates, and correspondingly change the resulting boundedness conclusion.
The key feature of our results is that they give non-trivial bounds for multilinear operators  whose integral kernels are measures which have Fourier transform in $L^{\infty}$ but satisfy no better uniform decay estimate. The positivity of the measures allows for the result to hold if the Fourier transform merely decays on the $d$-dimensional plane $\eta=-\xi$.  Before treating a more general class of non-translation invariant operators, we give proofs and counterexamples related to the translation invariant theorems and discuss  applications to continuous and discrete geometry.

The plan of the paper is as follows. In Sec. 2, we prove Theorem \ref{translationinvariant} and show that, despite its appearance, it is in some sense a bilinear theorem. In Sec. 3, we describe the a general framework of variations on the Falconer distance problem for $(k+1)$-point configurations by means of what we call \emph{$\Phi$-configurations} of points in $E\subset\mathbb R^d$, which include both $k$-simplices and their volumes. Applications of these results to Erd\"os-type problems in discrete geometry are given in Sec. 4. An extension of our main theorem to a nontranslation invariant setting and results about adjoints of multilinear operators are in Secs. 5 and 6, while Sec. 7 gives  a version of the regular value theorem for sets of fractional dimension. 

We would like to thank an anonymous referee for recommending expository improvements.


\vskip.25in

\section{Translation invariant proof and its intrinsic bilinearity}

\vskip.125in

\subsection{Proof of Theorem \ref{translationinvariant}}

We  assume that $f_j\ge 0$ are Schwartz functions. Using $f_j\ge 0$ in the first line and  the assumption \eqref{MainDecay} for the first inequality, one sees that
\begin{eqnarray*}
{||T_{\mu}(f_1, \dots, f_k)||}_{L^1(\mathbb{R}^d)}
&\leq & \prod_{j=3}^k {||f_j||}_{\infty}  \int \dots \int f_1(x-u^1) f_2(x-u^2) dx d\mu(u^1, \dots, u^k) 
\\
& =&\prod_{j=3}^k {||f_j||}_{\infty}   \int \dots \int f_1(y) f_2(y+u^1-u^2) dy d\mu(u^1, \dots, u^k)  
\\
&=& \prod_{j=3}^k {||f_j||}_{\infty}   \int \dots \int f_1(y) \int \widehat{f}_2(\xi)e^{2\pi i \xi\cdot(y+u^1-u^2)}d\xi dy d\mu(u^1, \dots, u^k)
\\
&=&
\prod_{j=3}^k {||f_j||}_{\infty} \int \widehat{f}_1(-\xi) \widehat{f}_2(\xi) \widehat{\mu}(-\xi, \xi, 0, \dots, 0) d\xi 
\\
& \lesssim & \prod_{j=3}^k {||f_j||}_{\infty} \int \left|\widehat{f}_1(-\xi)\right| \left|\widehat{f}_2(\xi)\right| (1+|\xi|)^{-\gamma} d\xi 
\\
& \lesssim &\prod_{j=3}^k {||f_j||}_{\infty} \left( \int {|\widehat{f}_1(\xi)|}^2 {(1+|\xi|)}^{-2\gamma_1} d\xi \right) \left( \int {|\widehat{f}_2(\xi)|}^2 {(1+|\xi|)}^{-2\gamma_2} d\xi \right).
\end{eqnarray*}

\subsection{Bilinearity of the multilinear estimates} 
Theorem \ref{translationinvariant} is inherently bilinear in nature: While the result is stated and used for multilinear operators, the assumption is distinctly bilinear in the sense that the number of derivatives gained over the trivial Holder estimate is based on the decay of the bilinear multiplier. We shall now see that this essentially unavoidable in the sense that a  translation invariant $k$-linear operator, $k \ge 3$, 
cannot in general gain the number of derivatives over the trivial Holder estimate corresponding to the {optimal} uniform decay of the $k$-linear multiplier. 
To see this, through a counterexample of geometric interest, let 
$$ T_d(f_1, \dots, f_d)(x)=\int \dots \int f_1(x-u^1) \dots f_d(x-u^d)\, d\Omega(u^1, \dots, u^d),$$ where, for any  $t\ne 0$, $d\Omega$ is the Leray measure on the smooth determinantal variety
$$\Sigma_t= \{(u^1, \dots, u^d) \in {\Bbb R}^d \times \dots \times {\Bbb R}^d: |\, u^1, \dots, u^d\, |:=det[u^1,\dots, u^d]=t \}. $$

\begin{theorem} \label{sharpnessexample} Let $T_k$ be as above and suppose that $d \ge 3$. Then 
\begin{equation} \label{dream} |\widehat{\Omega}(\xi^1, \dots, \xi^d)| \lesssim {(1+|\xi^1|+|\xi^2|+\dots+|\xi^d|)}^{-\frac{d^2-1}{2}}. \end{equation}
However, suppose that $p_j \ge 2$, $\frac{1}{p_1}+\dots+\frac{1}{p_d}=1$, $\gamma_j \ge 0$ and $\gamma_1+\dots+\gamma_d=\frac{d^2-1}{2}$. Then {\it no} estimate of the form 
\begin{equation} \label{leadtogold} {||T_d(f_1, \dots, f_d)||}_{L^1({\Bbb R}^d)} \lesssim {||f_1||}_{L^{p_1}_{-\gamma_1}({\Bbb R}^d)} \times \dots \times {||f_k||}_{L^{p_d}_{-\gamma_d}({\Bbb R}^d)} \end{equation} can hold, even on nonnegative functions. 
\end{theorem} 

To start the proof of   Theorem \ref{sharpnessexample}, we establish that the estimate (\ref{dream}) holds. Later we will  prove  that, if (\ref{leadtogold}) were to hold, that would imply that if the Hausdorff dimension of a subset of ${\Bbb R}^d$ is greater  than a number of the form  $d-1-\epsilon$, then the set of volumes determined by $(k+1)$-tuples of elements of $E$ is positive. This is absurd, since $E$ could be contained in a $(d-1)$-dimensional hyperplane.
\bigskip

 To see (\ref{dream}), consider $\Phi:\mathbb R^{d^2}\longrightarrow \mathbb R$, 
 $\Phi(x^1,\dots,x^d)=|\, x^1,\dots,x^d \, |$, so that $\Sigma_t=\Phi^{-1}(t)$. Then

\begin{equation}\label{Phi diff}
d\Phi(X^1,\dots,X^d)=\sum_{i=1}^d \big|x^1,\dots,x^{i-1},X^i,x^{i+1},\dots, x^d\big|.
\end{equation}

\noindent For a point $(x^1,\dots,x^d)\in\Sigma_t$, the vectors $x^1,\dots, x^d$ are linearly independent in $\R^d$; thus, $d\Phi(x)\ne 0$, since it is nonzero, e.g., when applied to vectors of the form $(X^1,\dots,X^d)=(0,\dots,0,c_ix^i,0,\dots,0)$, $c_i\ne 0$. Hence, $\Sigma_t$ is a smooth hypersurface. Note also that $SL(d,\R)$ acts transitively on $\Sigma_t$ through its diagonal action on $\R^{d^2}$,  since if $(x^1,\dots,x^d), \, (y^1,\dots, y^d)\in\Sigma_t$, then there exists a (unique) $A\in SL(\R,d)$ such that $Ay^i=x^i,\, 1\le i\le d$. Thus, to show that $\Sigma_t$ has nonzero Gaussian curvature everywhere, it suffices to consider its second fundamental form at a single point. Furthermore, since $\Sigma_t$ is a homothetic copy of $\Sigma_1$, we may assume that $t=1$. Thus, we may work at the  point $x_0=(e^1,\dots,e^d)\in\Sigma_1$, where the $e^i$ are the standard orthonormal basis for $\R^d$. 

Using the notation $\vec{X}=(X^1,\dots,X^d)\in T\R^{d^2}$ and $\big|\cdot,\cdot\big|_{ij}$ denoting the $(i,j)$-th $2\times 2$ minor of a $d\times 2$ matrix, differentiating  (\ref{Phi diff}) again yields that

\begin{eqnarray*}\label{Phi hessian}
d^2\Phi(\vec{X},\vec{X})&=& \sum_{1\le i<j\le d} \big| X^i,X^j\big|_{ij}, \\
&=& \sum_{1\le i<j\le d} X_i^iX_j^j-X_j^iX_i^j\\
&=& \sum_{1\le i<j\le d-1} X_i^iX_j^j-X_j^iX_i^j + \sum_{i=1}^d X_i^iX_d^d-X_d^iX_i^d.\\
\end{eqnarray*}
Now restrict this quadratic form to $T\Sigma_1$, using

\begin{equation}\label{tanspace}
T_{x_0}\Sigma_1=\Big\{(X^1,\dots,X^d) \, : \, X_1^1+\dots X_d^d=0\Big\}=\Big\{X_d^d=-\sum_{i=1}^{d-1} X_i^i\Big\};
\end{equation} this yields
\begin{equation}\label{Sec fund form}
\sum_{1\le i<j\le d-1} X_i^iX_j^j-\sum_{1\le i<j\le d-1}X_j^iX_i^j + \sum_{i=1}^{d-1} \Big(X_i^i\big(-\sum_{j=1}^{d-1}X_j^j\big)\Big) -\sum_{i=1}^{d-1}X_d^iX_i^d.
\end{equation}

From the second and fourth terms, we see that for $1\le i,j\le d,\, i\ne j$, the coordinate $X_i^j$ occurs (only) in the term $-X_j^iX_i^j$; each such term contributes 2 to the rank of the Hessian of $\Phi$. On the other hand, the $d-1$ coordinates  $X_i^i$, for $1\le i\le d-1$, occur as $-\sum_{1\le i\le j\le d-1} X_i^iX_j^j$, and  this quadratic form on $\R^{d-1}$ is represented by a $(d-1)\times (d-1)$ circulant matrix with first row $[1,\frac12,\dots,\frac12]$, which is easily seen to be nonsingular. The second fundamental form of $\Sigma_1$ at $x_0$ is thus nonsingular, showing that $\Sigma_1$ has non-zero Gaussian curvature there, and hence everywhere. Hence,   (\ref{dream}) follows
from the standard stationary phase estimate (see, e.g., \cite{St93}). 

We now prove that (\ref{leadtogold}) cannot hold, using Theorem \ref{applicationmama2} below.. 
If (\ref{leadtogold}) were to hold, it would imply that if the Hausdorff dimension of $E \subset {\Bbb R}^d$ is greater than 
$d-\frac{\gamma}{d}$, 
with $\gamma=\frac{d^2-1}{2}$, then the Lebesgue measure of ${\mathcal V}_d(E)$ has positive Lebesgue measure. This is not in general possible, since 
$$ d-\frac{\gamma}{d}=k-\frac{d^2-1}{2k}=\frac{d}{2}+\frac{1}{2d},$$ 
and this is smaller than $d-1$ when $d \ge 3$; since $E$ could be contained in a $(d-1)$-dimensional linear subspace, resulting in $\mathcal V_d(E)=\{0\}$, this is a contradiction. This proves Theorem \ref{sharpnessexample}.

\vskip.25in 


\section{Applications to problems in geometric measure theory} \label{sec geom}

\vskip.125in 

The classical Falconer distance problem, introduced in \cite{Fal86}, can be stated as follows: How large does the Hausdorff dimension of $E$ need to be to ensure that the Euclidean distance set $\Delta(E)=\{\, |x-y|\, : x,y \in E\}\subset\mathbb R$ has positive one-dimensional Lebesgue measure? This problem can be viewed as a continuous analogue of the Erd\H os distance problem (see \cite{Mat95} and the references  there). It is shown in \cite{Fal86}, using the set obtained by a suitable scaling of the thickened integer lattice that the best result we can hope for is the following. 
\begin{conjecture} \label{falconerconjecture} (Falconer distance conjecture) Let $E \subset {\Bbb R}^d$ with $dim_{{\mathcal H}}(E)>\frac{d}{2}$. Then the one-dimensional Lebesgue measure ${\mathcal L}^1(\Delta(E))>0$. 
\end{conjecture} 

The best  partial results known, due to Wolff \cite{W99} in the plane and to Erdo\u{g}an \cite{Erd05} in higher dimensions, say that ${\mathcal L}^1(\Delta(E))$ is indeed positive if the Hausdorff dimension $dim_{{\mathcal H}}(E)>\frac{d}{2}+\frac{1}{3}$. The proofs are Fourier analytic in nature and rely, at least in higher dimensions, on bilinear extension estimates. 

There are many possible multi-point configuration versions of the Falconer distance problem; perhaps the most immediate is the following. Let $E\subset\mathbb R^d,\, d\ge 2$, be compact; we may, for the sake of convenience, assume that $E \subset {[0,1]}^d$. For $2 \leq k \leq d+1$, call two $(k+1)$-tuples of vectors from the set $E \subset {\Bbb R}^d$ \emph{equivalent} if there exists a rigid motion that maps one $(k+1)$-tuple to the other, and let $T_k(E)$ denote the resulting set of equivalence classes. For $k=1$, $T_k(E)=\Delta(E)$, the distance set defined above, while for $k>1$, $T_k(E)$ can be thought of as the set of \emph{non-congruent $k$-simplices} in $E$. Note that, since rigid motions preserve distances, $T_k(E)$ can be naturally realized as  
$$\big\{(|x^i-x^j|)_{1\le i<j\le k+1} | x^1,\dots,x^{k+1}\in E\big\}\subset {\Bbb R}^{k+1 \choose 2},$$   
modulo permutations. Thus, the property of  having positive Lebesgue measure, ${\mathcal L}^{k+1 \choose 2}(T_k(E))>0$, is well-defined.  
It is then reasonable to ask whether one can find $0<s_0<d$ such that if $dim_{{\mathcal H}}(E)>s_0$,
then ${\mathcal L}^{k+1 \choose 2}(T_k(E))>0$. Relatedly, one can ask ``how likely" it is that  a given configuration can arise. 

More generally if $\Phi: \big( {\Bbb R}^d \big)^{k+1} \rightarrow {\Bbb R}^m$ is a smooth function, for some $ 1 \leq m \leq {d+1 \choose 2}$,
we define the \emph{set of $\Phi$-configurations of $E$} to be
$$ \Delta_{\Phi}(E):=\big\{\Phi(x^1, \dots, x^{k+1}): x^1,\dots,x^{k+1} \in E\big \}\subset\mathbb R^m.$$ 
Thus, for the noncongruent $k$-simplices above, $\Phi(x^1,\dots,x^{k+1})=\big(|x^i-x^j|\big)_{1\le i<j<k\le k+1}$.
We will focus on the \emph{translation invariant} case, relevant for configurations invariant under the additive structure of $\mathbb R^d$, for which $\Phi$ can be written as
$$ \Phi_0(x^1-x^{k+1}, \ldots , x^k-x^{k+1}) .$$
In analogy with the questions above, we ask how large the Hausdorff dimension of $E$ needs to be to ensure that 
$\mathcal L^m(\Delta_\Phi(E))>0$.
Letting $\nu$ be a Frostman measure supported on $E$, we also ask how likely it is for  $\Phi$-configurations to be near a fixed one  by asking how large the Hausdorff dimension of $E$ needs to be to ensure that uniform estimates of the form
\begin{equation} \label{mamaincidence}( \nu \times \dots \times \nu) \{(x^1, \dots, x^{k+1}): |\Phi(x^1, \dots, x^{k+1})-\vec{t}|<\epsilon \} \lesssim \epsilon^m \end{equation}
hold uniformly for $0<\epsilon<1$.

Our main result is the following. 

\begin{theorem} \label{applicationmama} Let $\Phi$ be translation invariant and $\Delta_{\Phi}(E)$ be defined as above. 
Let $t:=\vec{t}\in\R^m$ be a regular value of $\Phi$ and
$\mu_t$ be the Leray measure on the smooth surface
$$ \{(u^1, \dots, u^k): \Phi(u^1, \dots, u^k)=t \}.$$
Suppose that there exists a $\gamma > 0$ such that for all ordered $k$-tuples $\Xi_{j,\ell}$, $j\neq\ell$, with $-\xi$ in the $j$th entry and $\xi$ in the $\ell$th entry and $0$ in the remaining ones,
$$ |\widehat{\mu_t}(\Xi_{j,\ell})| \lesssim {(1+|\xi|)}^{-\gamma,}$$
and that for all ordered $k$-tuples $\Xi_{i}$ with $\xi$ in the $i$th entry and $0$ in the remaining entries,
$$|\widehat{\mu_t}(\Xi_{i})| \lesssim {(1+|\xi|)}^{-\gamma}.$$
Assume further that the Hausdorff dimension of $E$ is greater than $d-\frac{\gamma}{k}$. Then (\ref{mamaincidence}) holds and the $m$-dimensional Lebesgue measure  $\mathcal L^m(\Delta_{\Phi}(E))>0$. 
\end{theorem}

Theorem \ref{applicationmama} follows directly from Theorem \ref{translationinvariant}, Theorem \ref{adjoints}  and the following result,
which is proved in \S\S\ref{subsec 3.4}. 

\begin{theorem} \label{applicationmama2} Let $\Phi$ and $\mu_t$ be as in Thm. \ref{applicationmama}, and define
$$ T_{\Phi}(f_1, \ldots , f_k)(x) = \int\ldots\int f_1(x-u_1)\ldots f_k(x-u_k)\, d\mu_t(u_1,\ldots, u_k). $$
Suppose that  $p_j \ge 2$, $\gamma_j \ge 0$, $\gamma_1+\dots+\gamma_k=\gamma$, and 
\begin{equation} \label{partytime} {||T_{\Phi}(f_1, \dots, f_k)||}_{L^1({\Bbb R}^d)} \lesssim {||f_1||}_{L^{p_{\sigma(1)}}_{-\gamma_{\sigma(1)}}({\Bbb R}^d)} \dots 
{||f_k||}_{L^{p_{\sigma(k)}}_{-\gamma_{\sigma(k)}}({\Bbb R}^d)} \end{equation} for all permutations $\sigma$ and that the same estimates hold for each multilinear adjoint of $T_{\Phi}$. 
Then, if  $\dim_{\mathcal H}(E) > d-\frac{\gamma}{k}$,   (\ref{mamaincidence}) holds and hence  $\mathcal L^m(\Delta_{\Phi}(E))>0$. 

\end{theorem}

\vskip.125in 

\subsection{Distribution of simplices} Our main result, on the distribution of $k$-simplices in fractal sets, is the following; the proof will be given in \S\S\ref{subsec 3.5}. This result was previously established \cite{GI10} in the case of triangles in the plane, i.e., $d=k=2$. 

\begin{theorem} \label{kd} 
Let $d\ge 2$ and $1 \leq k \leq d$. Suppose that the Hausdorff dimension of a compact set $E \subset {\Bbb R}^d$ is greater than $s_0(k,d):=d-\frac{d-1}{2k}$. Let $\nu$ be  a Frostman measure supported on $E$, and ${\{t_{ij} \}}_{1 \leq i<j \leq k+1}$ a vector of positive real numbers. Then 
$$ (\nu \times \dots \times \nu) \{(x^1, \dots, x^{k+1}): t_{ij}-\epsilon \leq |x^i-x^j| \leq t_{ij}+\epsilon \} \lesssim \epsilon^{k+1 \choose 2}$$ with constant independent of $\epsilon$. 
Consequently, 
\begin{equation} \label{dimensionsimplex} {\mathcal L}^{k+1 \choose 2}(T_k(E))>0 \ \text{if} \ dim_{{\mathcal H}}(E)>d-\frac{d-1}{2k}. \end{equation} 
\end{theorem} 

\begin{remark} When $k=1$, Theorem \ref{kd} is known to be essentially sharp. See,e.g., \cite{IS10} and the references  there. When $k=d=2$, the estimate is also sharp, as was  shown in \cite{GI10}. In all the remaining cases, we believe that the Hausdorff exponents can be improved, and hope to address this in the future. \end{remark} 

\begin{remark} It is proved in \cite{EHI11} that ${\mathcal L}^{k+1 \choose 2}(T_k(E))>0$ if $dim_{{\mathcal H}}(E)>\frac{d+k+1}{2}$. One can check that the exponent in (\ref{dimensionsimplex}) is better when $d<k+2+\frac{1}{k-1}$. 
Moreover, the exponent $\frac{d+k+1}{2}$ is only $<d$ when $k<d-1$, whereas the exponent in (\ref{dimensionsimplex}) is always non-trivial.  \end{remark} 

\vskip.125in 

\subsection{Distribution of volumes of simplices} 

\vskip.125in 

In \cite{EHI11}, Erdo\u{g}an, Hart and the third listed author proved that if the Hausdorff dimension of $E \subset {[0,1]}^d$ is greater than $\frac{d+1}{2}$,  then the Lebesgue measure of the set of areas determined by triangles formed by pairs of points from $E$ and the origin is positive. In \cite{GIM11}, Mourgoglou, the second and the third listed authors proved an analogous result for volumes of simplices in ${\Bbb R}^3$ determined by three points of a given set and the origin. They prove that if $E\subset\mathbb R^3$ and
$dim_{{\mathcal H}}(E)>8/3$,
then the resulting set of volumes has positive Lebesgue measure in $\mathbb R$. 

In this section, we obtain a better exponent than the one in \cite{EHI11}  by considering $k$-dimensional volumes of simplices determined by $k+1$ points in $E$; a bilinear point of view once again plays a crucial rule. Our main result is the following. 

\begin{theorem} \label{volume} Define ${\mathcal V}_d(E)$ to be set of $d$-dimensional volumes determined by $(d+1)$-tuples of points from $E \subset {[0,1]}^d$, $d \ge 2$. Suppose that 
$dim_{{\mathcal H}}(E)$ is greater than $d-1+\frac{1}{2d}$ if $d$ is even, and greater than $d-1+\frac{1}{2(d-1)}$ if $d$ is odd. Then the Lebesgue measure of $\mathcal L^1\left({\mathcal V}_d\left(E\right)\right)>0$. 
\end{theorem} 

We shall prove this result in even dimensions and bootstrap it into odd dimensions using the following mechanism which is interesting in its own right. 

\begin{theorem} \label{bootstrap} Suppose that, whenever Hausdorff dimension of $E_{d-1} \subset {\Bbb R}^{d-1}$,  $d \ge 3$, is greater than $s_{d-1} \in (d-2,d-1)$, then the Lebesgue measure of ${\mathcal V}_{d-1}(E_{d-1})$ is positive. Then, if the Hausdorff dimension of $E_d \subset {\Bbb R}^d$ is greater than $s_d=s_{d-1}+1$, the Lebesgue measure of $\mathcal L^1\left({\mathcal V}_d\left(E_d\right)\right)>0$. 
\end{theorem} 

\begin{remark} Theorems \ref{volume} and \ref{bootstrap} will be proved in \S\S\ref{subsec 3.6}.
\end{remark}

\begin{remark} One immediate implication of Theorem \ref{bootstrap} is that if we could prove that $s_2=1$, then $s_d=d-1$, which would be the sharp exponent. Our current best exponent in two dimensions is $5/4$, which, in view of Theorem \ref{bootstrap} leads to the exponent $d-\frac{3}{4}$ in ${\Bbb R}^d$. We are able to obtain a better estimate in higher dimensions using multilinear operator bounds. Nevertheless, it would be reasonable to suppose that the ultimate resolution of the sharp exponent $d-1$ in $d$-dimensions will ultimately be accomplished by proving the sharp bound in dimension two and then using the boot-strapping mechanism of Theorem \ref{bootstrap}. 

\end{remark} 

\begin{remark} \label{kvolume} With a bit of work, our method extends to  $k$-dimensional volumes determined by $(k+1)$-tuples of vectors in ${\Bbb R}^d$. One can check that the main technical adjustment comes in the pigeon-holing argument at the end of the proof of Theorem \ref{volume}. In particular, the dimensional threshold exponent for $2$-dimensional volumes in ${\Bbb R}^d$ is 
$$ d-\frac{d-\frac{1}{2}}{2}=\frac{d}{2}+\frac{1}{4}.$$ 

The details are left to the interested reader. \end{remark} 

\vskip.125in 

\subsection{Distribution of angles} The following result is proved in \cite{IMP11} and is being included to illustrate the range of problems to which our method applies. 

\begin{definition} Let $E \subset {\mathbb R}^d$, $d \ge 2$. We say that an angle $\alpha \in [0, \pi]$ is {\it equitably represented} in $\mathcal{A}(E)$ if for every Frostman measure $\mu$ supported on $E$ and any $\epsilon>0$, 
\begin{equation}\label{def}
(\mu \times \mu \times \mu) \{(x^1,x^2,x^3): \alpha-\epsilon \leq \theta(x^1,x^2,x^3) \leq \alpha+\epsilon \} \lesssim  \epsilon,
\end{equation}
uniformly in $0<\epsilon<\\pi/2$, where $\theta:=\angle(x^2-x^1,x^3-x^1)$.
\end{definition}

\begin{theorem} \label{mainangle} Let $E$ be a compact subset of ${\mathbb R}^d$ of Hausdorff dimension greater than $\frac{d+1}{2}$. Then every $\alpha \in [0, \pi]$ is equitably represented in $\mathcal{A}(E)$. 
\end{theorem}

\begin{corollary} \label{mainanglecorollary}  Let $E$ be a compact subset of ${\mathbb R}^d$ of Hausdorff dimension greater than $\frac{d+1}{2}$. Then the Lebesgue measure of ${\mathcal A}(E)$ is positive. 
\end{corollary} 

Theorem \ref{mainangle} can be recovered from Theorem \ref{applicationmama} via the following lemma proved in \cite{IMP11}. 
\begin{lemma} \label{anglelemma} Let 
\begin{equation} \label{bendover2} \widehat{\mu}_t(\xi, \eta)=\iiint e^{-2 \pi i (u \cdot \xi+ \lambda \theta u \cdot \eta)} \psi(|u|) \psi_0(\lambda) d\Omega_{t, \frac{u}{|u|}}(\theta)du d\lambda, \end{equation} where $d\Omega_{t, \frac{u}{|u|}}(\theta)$ is the restriction of the Haar measure on $SO(d)$ to $\Omega_{t, \frac{u}{|u|}}$ and $\psi, \psi_0$ are smooth cut-off functions. Then 
$$|\widehat{\mu}_t(\xi, \eta)| \lesssim \frac{1}{(1+|\xi|+|\eta|)^{d-1}}. $$
\end{lemma} 

\vskip.125in 


\subsection{Proof of Theorem \ref{applicationmama2}} \label{subsec 3.4}

\vskip.125in

Let $\nu$ be a Frostman measure, supported on $E \subset {\mathbb R}^d$. Set $s=\dim_{\mathcal{H}}(E)$. Let $\rho$ be a non-negative, smooth function, equal to 1 on $[-\frac{1}{4}, \frac{1}{4}]$, supported in $[-1,1]$ with $\| \rho \|_{L^1(\mathbb{R}^d)} = 1$, and $\nu^{\delta} := \nu * \rho_{\delta}$, where $\rho_{\delta}(x) = \delta^{-d}\rho(\frac{x}{\delta})$, the resulting smooth approximation of $\nu$ as $\delta$ tends to $0$. We will establish the bound
$$ \frac{1}{\epsilon^m}\ (\nu^{\delta} \times \dots \times \nu^{\delta}) \{(x^1, \dots, x^{k+1}): |\Phi(x^{k+1} - x^1, \ldots , x^{k+1} - x^k)-\vec{t}\,|<\epsilon \} \lesssim 1,$$
uniformly in  $\delta$, and thus by passing to the limit we establish the theorem.
Write
$$ \frac{1}{\epsilon^m}\ (\nu^{\delta} \times \dots \times \nu^{\delta}) \{(x^1, \dots, x^{k+1}): |\Phi(x^{k+1} - x^1, \ldots , x^{k+1} - x^k)-\vec{t}\,|<\epsilon \}$$ 
$$=\epsilon^{-m}\ \int \ldots \int_{\{(x^1, \dots, x^{k+1}): |\Phi(x^{k+1} - x^1, \ldots , x^{k+1} - x^k)-\vec{t}|<\epsilon \}} d\nu^{\delta}(x_1)\ldots d\nu^{\delta}(x_k) d\nu^{\delta}(x_{k+1})$$
$$=\epsilon^{-m}\ \int \ldots \int_{\{(u^1, \dots, u^{k+1}): |\Phi(u^1, \ldots ,  u^k)-\vec{t}|<\epsilon \}} d\nu^{\delta}(u_1)\ldots d\nu^{\delta}(u_k) d\nu^{\delta}(x_{k+1})$$
\begin{equation} \label{fp}=\langle T_{\Phi}(\nu^{\delta}, \ldots, \nu^{\delta}), \nu^{\delta} \rangle ,\end{equation}
and $\langle \cdot,\cdot \rangle$ is the $L^2(\mathbb{R}^d)$ inner product.

Now define
$$ F(\alpha_1, \ldots, \alpha_k, \alpha_{k+1}) := \langle T_{\Phi}(\nu^{\delta}_{\alpha_1},\ldots , \nu^{\delta}_{\alpha_k}), \nu^{\delta}_{\alpha_{k+1}} \rangle $$
where, initially defined for $\text{Re}(\beta) > 0$, 
\begin{equation}\label{spongebob}
\nu^{\delta}_{\beta}(x) := \frac{2^{\frac{d-\beta}{2}}}{\Gamma\left(\beta/2 \right)}(\nu^{\delta} * |\cdot |^{-d+\beta})(x)
\end{equation}
is extended to the complex plane by analytic continuation. Since $\nu^{\delta}_{\beta}$ is smooth and we are in a compact setting, we have trivial bounds on $F(\alpha_1, \dots, \alpha_{k+1})$ with constants depending on $\delta$. Observe that $\widehat{\nu^{\delta}_{\beta}} (\xi) = C_{\beta,d}\widehat{\nu}(\xi)\widehat{\rho}(\delta\xi)|\xi|^{-\beta}$ where
\begin{equation}\label{blowup}
C_{\beta,d} = \frac{ \pi^{\frac{d}{2}} }{\Gamma(\frac{d-\beta}{2})}.
\end{equation}

See, e.g., \cite[p. \!192]{GS58} for this and related calculations. By Plancherel, $\nu_{\beta}^\delta$ is an $L^2({\mathbb R}^d)$ function with bounds depending on $\delta$. Taking the modulus in (\ref{spongebob}), we see that
$$|\mu_\beta^\delta(x)|\leq \left| \frac{ 2 ^{\frac{d-\beta}{2}}}{\Gamma\left(\frac{\beta}{2} \right)} \right| (\nu^{\delta} * |\cdot |^{-d+\text{Re}(\beta)})(x) = \frac{\Gamma(\text{Re}(\frac{\beta}{2}))}{|\Gamma(\frac{\beta}{2})|}\nu_{\text{Re}(\beta)}^\delta(x)$$
and note, using (\ref{spongebob}), that the right hand side is non-negative. To be more clear then the inequality clearly holds for $\text{Re}(\beta) > 0$ but can then be extended to $\text{Re}(\beta) \leq 0$ using the regularization needed for the analytic extension, see pages 50-56 and pages 71-72 in \cite{GS58}.
Theorem \ref{applicationmama2} would thus follow if we could show that, whenever 
$dim_{{\mathcal H}}(E)>d-\frac{\gamma}{k}$, then $F(0, \dots, 0) \lesssim 1$.

Instead of pursuing the bound $F(0, \dots, 0) \lesssim 1$ directly we will obtain it through interpolation. The key tool is the following straightforward multilinear generalization of the three lines lemma.

\begin{lemma}\label{multi3linelemma}
Let $a_1,\ldots,a_k\in\mathbb{R}^{n}$ and $f(z_1,\ldots,z_n)$ be a bounded function of $z_j = x_j + i y_j$, $j=1,\ldots,n$, defined on the set
$$\lbrace (x_1 + i y_1,\ldots, x_n + i y_n) : (x_1,\ldots,x_n) = t_1 a_1 + \ldots + t_k a_k, 0\leq t_1,\ldots,t_k \leq 1, t_1 + \ldots + t_k = 1 \rbrace ,$$
holomorphic in the interior and continuous on the whole set. If
$$M(x_1,\ldots,x_n) = \sup\limits_{y_1,\ldots,y_n}|f(x_1+i y_1,\ldots, x_n + i y_n)| $$
then for any $(x_1,\ldots,x_n)=t_1 a_1 + \ldots + t_k a_k$ where $0\leq t_1,\ldots,t_k \leq 1$ and $t_1 + \ldots + t_k = 1$ we have
$$ M(x_1,\ldots,x_n) \leq M(a_1)^{t_1}\cdot\ldots\cdot M(a_k)^{t_k} .$$
\end{lemma}

We shall bound $F(\alpha_1, \dots, \alpha_{k+1})$, with a constant independent of $\delta$, for every possible $(k+1)$-tuple of complex numbers where
\begin{equation} \label{alphadog} \text{Re}(\alpha_i) = \frac{\gamma}{kp_{i}'}-\gamma_i\end{equation}
for $1\leq i \leq k$ and
\begin{equation} \label{alphadoginf} \text{Re}(\alpha_{k+1}) = \frac{\gamma}{k}\end{equation}
and all permutations of such numbers. By Lemma \ref{multi3linelemma} we then obtain a bound for $F(\alpha_1, \dots, \alpha_{k+1})$, independent of $\delta$, in the convex hull of these numbers. It is not difficult to see that the origin in ${\mathbb C}^{k+1}$ is contained in the convex hull as the sum of the real parts of these numbers is $0$ and, for each coordinate, all the $\alpha_i$'s appear equally often.

We get boundedness for each permutation of the $\alpha$'s from precisely one assumption in the statement of the theorem. Without loss of generality, $\alpha_i$ equals the expression in (\ref{alphadog}) for $1\leq i \leq k$ and $\alpha_{k+1}=\frac{\gamma}{k}$. We have 
\begin{equation} \label{battlebegins} F(\alpha_1, \dots, \alpha_{k+1})= \big\langle \,T(\nu_{\alpha_1}^{\delta}, \nu_{\alpha_2}^{\delta}, \dots, \nu_{\alpha_k}^{\delta}), \nu_{\alpha_{k+1}^{\delta}}\big\rangle.\end{equation} 

We need the following simple observation. 
\begin{lemma} \label{Linfinity} If $\nu$ is a Frostman measure on a set  $E\subset\mathbb R^d$ with $dim_{{\mathcal H}}(E)> s$, then 
$$ {||\nu_{\alpha}^{\delta}||}_{\infty} \lesssim 1 \ \text{if} \ Re(\alpha)=d-s.$$ 
\end{lemma} 

To prove the lemma,  observe that if $Re(\alpha)=d-s$, 

\begin{equation}\nonumber
|\nu_{\alpha}^{\delta}(x)| \leq \int {|x-y|}^{-s} d\nu^{\delta}(y) \approx \sum_j 2^{js} \int_{|x-y| \approx 2^{-j}} d\nu^{\delta}(y) \lesssim  \sum_j 2^{js} 2^{-j \cdot dim_{{\mathcal H}}(E)},
\end{equation} 
and this  is $\lessapprox 1$, since $\nu$ is a Frostman measure on $E$ and  $dim_{{\mathcal H}}(E)>s$. 

We now bound
\begin{equation} \label{fpp} 
|F(\alpha)| \leq {||T(\nu_{\alpha_1}^{\delta}, \nu_{\alpha_2}^{\delta}, \dots, \nu_{\alpha_k}^{\delta})||}_{1} {||\nu_{\alpha_{k+1}}^{\delta}||}_{\infty} \end{equation} \
$$ \lesssim {||T(\nu_{\text{Re}(\alpha_1)}^{\delta}, \nu_{\text{Re}(\alpha_2)}^{\delta}, \dots, \nu_{\text{Re}(\alpha_k)}^{\delta})||}_{1} {||\nu_{\text{Re}(\alpha_{k+1})}^{\delta}||}_{\infty},$$ 
where the implicit constant depends on terms with gamma functions. Note that if, e.g., $\frac{\gamma}{2k}\in\mathbb Z$, then our bounds blow up. This is however not an obstacle, because we can instead consider $\frac{\gamma}{2k}+\epsilon$, where $\epsilon>0$ is chosen so small that   we still have that $\dim_{\mathcal H}(E)>d-\frac{\gamma}{k}+\epsilon$.
\medskip

By Lemma \ref{Linfinity}, the expression in (\ref{fpp}) is bounded by 
\begin{equation} \label{mfp} {||T(\nu_{\text{Re}(\alpha_1)}^{\delta}, \nu_{\text{Re}(\alpha_2)}^{\delta}, \dots, \nu_{\text{Re}(\alpha_k)}^{\delta})||}_1. \end{equation} 
By assumption this expression is bounded by
\begin{equation}\label{e} \prod_{j=1}^k {||\nu_{\text{Re}(\alpha_j)}^{\delta}||}_{L^{p_j}_{-\gamma_{j}}({\Bbb R}^d)} \end{equation} and all we need to do is to establish the boundedness of 
$||\nu_{\text{Re}(\alpha_j)}^{\delta}||_{L^{p_j}_{-\gamma_{j}}({\Bbb R}^d)}$ for each $j$. Note that we can rewrite this as $$||\nu_{\text{Re}(\alpha_j)+\gamma_j}^{\delta}||_{L^{p_j}({\Bbb R}^d)} = ||\nu_{\frac{\gamma}{k}\left(1 - \frac{1}{p_j}\right)}^{\delta}||_{L^{p_j}({\Bbb R}^d)}.$$
We now calculate using H\"{o}lder
\begin{align*}
||\nu_{\frac{\gamma}{k}\left(1 - \frac{1}{p_j}\right)}^{\delta}||_{L^{p_j}({\Bbb R}^d)}^{p_j} &\lesssim \int\left(\nu^{\delta}(y) |x-y|^{-d+\frac{\gamma}{k}(1-\frac{1}{p_j})} \right)^{p_j}dx \\
&= \int\left(\int \left( \nu^{\delta}(y) |x-y|^{-d+\frac{\gamma}{k}} \right)^{1-\frac{2}{p_j}} \left(\nu^{\delta}(y) |x-y|^{-d+\frac{\gamma}{2k}} \right)^{\frac{2}{p_j}} dy\right)^{p_j}dx \\
&\lesssim \int \left( \int \nu^{\delta}(y) |x-y|^{-d+\frac{\gamma}{k}} dy \right)^{p_j - 2} \left( \int \nu^{\delta}(y) |x-y|^{-d+\frac{\gamma}{2k}} dy \right)^{2} dx \\
&\lesssim ||\nu_{\frac{\gamma}{k}}^{\delta}||_{L^{\infty}({\Bbb R}^d)}^{p_j - 2} ||\nu_{\frac{\gamma}{2k}}^{\delta}||_{L^{2}({\Bbb R}^d)}^2
\end{align*}
and this is bounded using Lemma \ref{Linfinity} and the following lemma.

\begin{lemma} \label{Ltwo} If $\nu$ is a Frostman measure on a set  $E\subset\mathbb R^d$ with $dim_{{\mathcal H}}(E)> s$, then
$$ {||\nu_{\alpha}^{\delta}||}_{L^2(\mathbb{R}^d)} \lesssim 1 \ \text{if} \ Re(\alpha)=\frac{s}{2}.$$ 
\end{lemma} 

To prove the lemma,  first recall the energy integral,
$$ I_s(\nu)=\int \int {|x-y|}^{-s} d\nu(x) d\nu(y).$$
A standard calculation shows
$$ {||\nu_{\alpha}^{\delta}||}_{L^2(\mathbb{R}^d)} \lesssim \left(I_{d-2\text{Re}(\alpha)}(\nu^{\delta})\right)^{1/2} =   \left(I_{d-s}(\nu^{\delta})\right)^{1/2} $$
and since $\nu$ is supported on a set of Hausdorff dimension greater than $s$, we can bound
$ I_{d-s}(\nu^{\delta}) \lesssim 1 $,
with a bound independent of $\delta$.

\vskip.125in 

\subsection{Proof of Theorem \ref{kd}} \label{subsec 3.5}

\vskip.125in

In light of Theorem \ref{applicationmama} it suffices to establish the following multiplier estimate:

\begin{lemma} \label{iscool} On ${\Bbb R}^d \times {\Bbb R}^d$, $d\ge 3$,  let 
$$K(x,y)=\delta(|x|-1)\delta(|y|-1)\delta(|x-y|-1)\simeq\delta(|x|^2-1)\delta(|y|^2-1)\delta(x\cdot y-\frac12).$$

Then 
$$|\widehat{K}(\xi,\eta)|\lesssim (1+|\xi|+|\eta|)^{-\frac{d-1}2}.$$ 
\end{lemma} 

To see why this lemma is sufficient then note that according to Theorem \ref{applicationmama} we need to consider
$$ \int\limits_{\substack{|u^i| = t_i \\|u^i - u^j| = t_{ij}}}e^{2\pi i (u^1\cdot\xi^1 + \ldots + u^k\cdot \xi^k)} du^1 \ldots du^k $$
where we are using that this point configuration is translation invariant. Further note that according to Theorem \ref{applicationmama} we only need to understand the decay of this integral for very particular choices of the Fourier variables, to be precise, we only need to consider $(\xi^1,\ldots,\xi^k)$ that are of the type $(\xi,-\xi,0,\ldots,0)$ or $(\xi,0,\ldots,0)$ and permutations thereof. In order to tackle both cases in one go we consider a generic setup where we have two Fourier variables $\xi$ and $\eta$ and the remaining $k-2$ Fourier variables are all $0$. For the space variables associated to the Fourier variables that are $0$ we simply integrate through. Note that in each case we get that the variable we are integrating through lies on one more sphere than the number of variables are left and the distances between the centers of these spheres are fixed and depend only on $t_i$'s and $t_{ij}$'s. Thus we get a constant contribution that depends on the $t_i$'s and $t_{ij}$'s, the dimension $d$ and $k\leq d$. It is worth recalling that all entries in the configuration vector $\vec{t}$ are positive so we face no degenerate situations. At the very end the constant we obtain this way can be completely dominated by a constant depending only on the dimension $d$ and the diameter of the compact set $E$ that we began with. This happens because spheres with diameter larger than that of the compact set $E$ that we began with are guaranteed to have no contribution. Finally for the remaining two variables it is clear that it is sufficient to show the boundedness of an expression as in Lemma \ref{iscool}, with the the sole exception that the corresponding $t_i$'s and $t_{ij}$'s might not be equal to $1$. However, it is clear from the proof of the lemma that we get a similar estimate in the general case where we just introduce an extra constant that only depends on the diameter of $E$ and the dimension $d$. 

To prove the lemma, we use a partition of unity to decompose $K$ into a finite sum of terms supported on products of small spherical caps about points $x^j,y^j,\, 1\le j\le N_d$. By rotation invariance of $K$, we can assume that the basepoints are $x^0=(0,\zero',1)$ and $y^0=(\frac{\sqrt3}2,\zero',\frac12)$, where we write $\R^d\owns x=(x_1,x',x_d)$. Introduce local coordinates on $\sd$ near $x^0,y^0$, resp.,
\bes
x(u)&=&(u,1-\frac{|u|^2}2)+\O(|u|^3),\, u=(u_1,u')\in\R^{d-1},\, |u|\le \epsilon,\hbox{ and }\\
y(v)&=&(\frac{\sqrt3}2 +v_1,v',\frac12-\st v_1-4v_1^2-|v'|^2)+\O(|v|^3),\, v=(v_1,v')\in\R^{d-1},\, |v|<\epsilon.
\ees
All calculations that follow will be modulo $\O^3:=\O(|u,v|^3)$. 

The measure $K^0$ is a smooth multiple of surface measure on the $(2d-3)$-dimensional manifold $\{(u,v): x(u)\cdot y(v)=\frac12\}\subset \R^{d-1}\times\R^{d-1}$, pushed forward under the parametrization map $(u,v)\rightarrow (x(u),y(v))$. However, 
\bes
x(u)\cdot y(v)-\frac12&=&\frac{\st}2 u_1-\st v_1+u\cdot v -\frac{|u|^2}4-4v_1^2-|v'|^2 +\O^3=0\Leftrightarrow \\
(\st-u_1)v_1 + 4v_1^2&=&\frac{\st}2 u_1+u'\cdot v'-\frac{|u|^2}4-|v'|^2+\O^3.
\ees
We use the  quadratic terms in the implicit function theorem in one variable,
\be
a_1 s+ a_2 s^2=t\implies s=a_1^{-1}t-a_1^{-3}a_2t^2+\O(t^3),\, s,t\searrow 0,
\ee
to solve for $v_1$ in terms of $u_1$, with $u',v'$ as parameters:
\be
v_1=\frac12 u_1-\frac{13\st}{36} u_1^2-\frac{\st}{12} |u'|^2-\frac1\st |v'|^2+\frac{1}{\st}u'\cdot v' +\O^3,
\ee
so that also $v_1^2=\frac14 u_1^2+\O^3$. Hence, the contribution $dK^0$ of this pair of spherical caps  to $dK$ is the pushforward of $du_1\, du'\, dv'$ under
\bes \big(u_1,u',v')&\rightarrow& (u_1,u',1-\frac12 u_1^2-\frac12 |u'|^2;\\
& & \frac{\st}2+\frac12u_1-\frac{13\st}{36} u_1^2-\frac{\st}{12}|u'|^2-\frac1{\st}u'\cdot v',\\
& &\quad v', \frac12-\frac\st2u_1+\frac76u_1^2+\frac14|u'|^2-u'\cdot v'\big).
\ees
Thus, for $\chi$ a suitable smooth cutoff supported near $(\zero,\zero')$, 
\be\label{eqn osc int}
\widehat{K^0}(\xi,\eta)=\int e^{-i\Psi(u_1,u',v')} \chi(u_1,u',v') du_1\, du'\, dv',
\ee
where the phase function is
\bes
\Psi&:=&\xi_1u_1+\xi'\cdot u' +\xi_d\big(1-\frac{|u|^2}2\big)+\eta'\cdot v'\\
& &+\eta_1\big(\frac{\st}2+\frac12u_1-\frac{13\st}{36} u_1^2-\frac{\st}{12}|u'|^2-\frac1\st |v'|^2+\frac1\st u'\cdot v'\big)\\
& &+\eta_d\big(\frac12-\frac\st2u_1+\frac76u_1^2+\frac14|u'|^2-u'\cdot v'\big),
\ees
and the (linear) dependence of $\Psi$ on $(\xi,\eta)$ is notationally suppressed. For $(\xi,\eta)$ fixed, $\Psi$ has a unique critical point; conversely, if we fix $(u,v')$, then the system of linear equations

\be
d_{u_1}\Psi(u,v')=0,\, d_{u'}\Psi(u,v')=d_{v'}\Psi(u,v')=0
\ee
 in $(\xi,\eta)$ is of maximal rank, with a 3-dimensional solution space, corresponding to the conormal bundle of $\Sigma$ at 
 $\left(x\left(u\right),y\left(v_1\left(u,v'\right),v'\right)\right)$. Without loss of generality, we consider $(u_1,u',v')=(0,\zero',\zero')$
 and show the claimed decay rate for $\widehat{K^0}$ in the conormal plane there, $N^*\Sigma=\{\xi_1+\eta_1-\frac\st2\eta_d=0,\, \xi'=\eta'=0\}$.
 Splitting off $u_1$ and then pairing $u_i,v_i,\, 2\le i\le d-1$, one computes that the Hessian of $\Psi$ at $(0,\zero',\zero')$ is 
 \be
 p(\xi,\eta)\I_{1\times1}\oplus \bigoplus_2^{d-1} \Big[
 \begin{matrix}
 q_{11}(\xi,\eta) & q_{12}(\xi,\eta)\\
 q_{21}(\xi,\eta) & q_{22}(\xi,\eta)
 \end{matrix}
 \Big],
 \ee
where 
\be
p(\xi,\eta)=-\xi_d-\frac{13\st}{18}\eta_1+\frac73\eta_d
\ee
and
\be
q_{11}=-(\xi_d+\frac\st6\eta_1-\frac12\eta_d),\, q_{12}=q_{21}=\frac1\st\eta_1-\eta_d,\,  q_{22}=-\frac2\st\eta_1.
\ee
This has determinant $\pm p(\xi,\eta)\cdot q(\xi,\eta)^{d-2}$, where 
\be\label{eqn hess}
q(\xi,\eta)=
 \Big|
 \begin{matrix}
 q_{11}(\xi,\eta) & q_{12}(\xi,\eta)\\
 q_{21}(\xi,\eta) & q_{22}(\xi,\eta)
 \end{matrix}
 \Big|.
\ee
If we restrict to the hyperplane   $\Pi:=\{p=0\}$ by setting $\xi_d=-\frac{13\st}{18}\eta_1+\frac73\eta_d$, then 
\be\label{eqn q}
q=-\frac{13}9\eta_1^2+\frac{17\st}{9}\eta_1\eta_d-\eta_d^2,
\ee
which has negative discriminant and thus is a negative-definite quadratic form. Hence, within $N^*\Sigma$,  the cone $\Gamma:=\{q(\xi,\eta)=0\}$ intersects $\Pi$ only at the origin. Off of $\Gamma$,
$\hbox{rank}\, \Psi''\ge 2(d-2)$, so that stationary phase gives $|\widehat{K^0}(\xi,\eta)|\lesssim (1+|\xi|+|\eta|)^{-(d-2)}$, and this estimate is uniform if $dist((\xi,\eta),\Gamma)\ge c(|\xi|+|\eta|)$; for $d\ge 3$, this is at least as good as in the statement of the proposition. On the other hand, by (\ref{eqn hess}), (\ref{eqn q}) and the comment following it, the rank  of the Hessian at points of $\Pi$ is at least $1+(d-2)=d-1$, yielding $|\widehat{K^0}(\xi,\eta)|\lesssim (1+|\xi|+|\eta|)^{-\frac{d-1}2}$, and this is uniform on a conic neighborhood of $\Pi$.
By standard proofs of stationary phase, these estimates are uniform both in the conormal directions above the base point, and also as the base point is varied, completing the proof of Lemma \ref{iscool}.

\vskip.25in

\subsection{Proofs of Thms. \ref{volume} and  \ref{bootstrap}} \label{subsec 3.6}

\vskip.125in 

We first prove Theorem \ref{bootstrap}. If $E \subset {\Bbb R}^d$, $d \ge 3$,  and $dim_{\mathcal H}(E)>s_d$, then there exists a $(d-1)$-plane $H$ such that $dim_{\mathcal H}(E \cap H)>s_d-1$; see, e.g., \cite{Mat95}. Define $s_{d-1}=s_d-1$. By assumption, ${\mathcal V}_{d-1}(E \cap H)>0$. Let $H'$ denote a $(d-1)$-plane parallel to $H$, but not equal to $H$, containing at least one point of $E$, denoted by $z$. Such a plane must exist since $s_d>d-1$ by assumption. For every $d$-tuple in $E \cap H$ that contributes a non-zero element to ${\mathcal V}_{d-1}(E \cap H)$, form a $(d+1)$-tuple in ${\Bbb R}^d$ by adjoining $z$. The volumes of the resulting $(d+1)$-tuples are all distinct and form a set of positive one-dimensional Lebesgue measure. This completes the proof of Theorem \ref{bootstrap}.

To prove Theorem \ref{volume} we first observe that Theorem \ref{applicationmama} goes over with no substantive changes if instead of assuming that each of vectors in the $(k+1)$-point configuration under consideration being contained in a single set $E$, we have $k+1$ sets and $x^j \in E_j$, $j=1, 2, \dots, k+1$. The dimensional assumption then becomes $dim_{{\mathcal H}}(E_j)>d-\frac{\gamma}{k}$ for each $1 \leq j \leq k+1$.

Now, in view of Theorem \ref{applicationmama}, using the translation invariant configuration function \linebreak$\Phi=det(x^1-x^{d+1},\dots, x^d-x^{d+1})$, the even dimensional case of Theorem \ref{volume} would follow from:
\begin{lemma} \label{partialdeterminant} Let 
\[ f(x,y)=\det \left( \begin{array}{ccc}
x \\
y \\
z^1 \\
z^2 \\
  .   \\
  .  \\
  . \\
  
z^{d-2} \end{array} \right), \] where $z^j \in E^j$ as above.
Then in any compact $B \subset {\Bbb R}^{2d}$ and $t \not=0$, the hypersurface
$ \{(x,y): f(x,y)=t \}$ has at least $2d-1$ non-vanishing principal curvatures, with bounds independent of $z^j$, and thus 
the Fourier transform of the Leray measure on this surface decays uniformly of order $-\frac{2d-1}{2}$ at infinity. 
\end{lemma} 

To prove Lemma \ref{partialdeterminant}, consider the submatrix \[ \left( \begin{array}{ccc}
x_1 \ x_2 \\
y_1 \ y_2 \\ 
\end{array} \right), \] 
When computing the determinant of the whole matrix, we get $x_1y_2-x_2y_1$ times the determinant of the matrix obtained by covering up the first two columns and the first two rows. The determinant of this matrix is non-zero by definition. Thus 
$$ f(x,y)=c_{12}(x_1y_2-x_2y_1)+c_{34}(x_3y_4-x_4y_3)+\dots+c_{d-1 d} (x_{d-1}y_d-x_dy_{d-1}),$$ where $c_{ij}$s are non-zero. A simple rotation $(x,y)\to (u,v)$ transforms   $f(x,y)$ into 
$$ c_{12}(u_1^2-v_2^2-u_2^2+v_1^2)+\dots+c_{d-1 d}(u_{d-1}^2-v_d^2-u_d^2+v_{d-1}^2),$$ a non-degenerate form, so the number of non-vanishing principal curvatures is indeed $2d-1$. The conclusion of Lemma \ref{partialdeterminant} follows. 

This takes care of the even dimensional case of Theorem \ref{volume}; the odd dimensional case follows at once by combining this with Theorem \ref{bootstrap}.

\vskip.25in



\section{Applications to discrete geometry} 
\label{discreteapplications}

\vskip.125in 

The purpose of this section is to use a variant of a ``continuous-to-discrete" mechanism, developed in \cite{HI05,IL05,I08,IRU10}, to show that Thms. \ref{kd} and  \ref{volume} 
yield interesting discrete analogues. We work in the setting of \emph{$s$-adaptable} sets, defined below. 
These are more general than the so-called {\it homogeneous} sets, used widely in geometric combinatorics, e.g., \cite{SV04} and the references there. 

\begin{definition} (\cite{IRU10}) Let $P$ be a set of $n$ points contained in ${[0,1]}^d$, $d \ge 2$. Let $\chi_{B^{\, p}_{n^{-\frac{1}{s}}}}(x)$ be the characteristic function of the ball of radius $n^{-\frac{1}{s}}$ centered at $p$, and define the measure
\begin{equation} \label{pizdatayamera} d \mu^s_P(x)=n^{-1} \cdot n^{\frac{d}{s}} \cdot \sum_{p \in P} \chi_{B^{\, p}_{n^{-\frac{1}{s}}}}(x)\, dx. \end{equation} 
We say that $P$ is \emph{$s$-adaptable}  if 

$$ I_s(\mu^s_P)=\int \int {|x-y|}^{-s} d\mu^s_P(x) d\mu^s_P(y)<\infty.$$ 

\end{definition} 

\vskip.125in

This is equivalent to the statement 

\begin{equation} \label{discreteenergy} n^{-2} \sum_{p \not=p' \in P} {|p-p'|}^{-s} \lesssim 1.\end{equation} 

\vskip.125in 

To understand this condition in  clearer geometric terms, suppose that $P$ comes from a \linebreak$1$-separated set $A\subset\mathbb R^d$, rescaled down by its diameter so as to be contained in $[0,1]^d$. Then the condition (\ref{discreteenergy}) takes the form 

\begin{equation} \label{discreteenergylarge} n^{-2} \sum_{a \not=a' \in A} {|a-a'|}^{-s} \lesssim {(diameter(A))}^{-s}. \end{equation} 
This says $P$ is $s$-adaptable if it is a scaled $1$-separated set where the expected value of the distance between two points raised to the power $-s$ is comparable to the value of the diameter raised to the power of $-s$. Thus, for a set to be $s$-adaptable, clustering is not allowed to be too severe, on average. 

To put it in more technical terms, $s$-adaptability means that a discrete point set $P$ can be thickened into a set which is uniformly $s$-dimensional in the sense that its energy integral of order $s$ is finite. Unfortunately, it is shown in \cite{IRU10} that there exist finite point sets which are not $s$-adaptable for certain ranges of the parameter $s$. The point is that the notion of Hausdorff dimension is much more subtle than the simple ``size" estimate. However, many natural classes of sets are $s$-adaptable. For example, all homogeneous sets, studied by Solymosi and Vu \cite{SV04} and others, are $s$-adaptable for any $0<s<d$. See also \cite{IJL09} where $s$-adaptability of homogeneous sets is used to extract discrete incidence theorems from Fourier type bounds. 

While $s$-adaptability is a restriction, we will see below that combining this notion with analytic methods allows one to prove robust theorems, involving \emph{neighborhoods} of a class of geometric objects, something that is not typically possible using combinatorial methods. 

\subsection{Simplices determined by discrete sets} Before we state the discrete result that follows from Theorem \ref{kd}, let us briefly review what is known. If $P$ is set of $n$ points in ${[0,1]}^2$, let $u_{2,2}(n)$ denote the number of times a fixed triangle can arise among points of $P$. It is not hard to see that 
\begin{equation} \label{st} u_{2,2}(n)=O(n^{\frac{4}{3}}).\end{equation} 
This follows  from the fact that a single distance cannot arise more than $O(n^{\frac{4}{3}})$ times, which, in turn, follows from the celebrated Szemeredi-Trotter incidence theorem; see \cite{BMP05} and the references  there. By the pigeon-hole principle, one can conclude that 

\begin{equation} \label{distance} \# T_2(P) \gtrsim \frac{n^3}{n^{\frac{4}{3}}}=n^{\frac{5}{3}}.\end{equation}
(Recall that $T_k(P)$ is the set of $k$-point configurations of $P$, as in Sec. \ref{sec geom}.)

However, one can do quite a bit better as far as the lower bound on $\# T_2(P)$ is concerned. It is shown in \cite[p. 263]{BMP05} that 
$$ \# T_2(P) \gtrsim n \cdot \# \big\{|x-y|: x,y \in P\big\}.$$ 

Guth and Katz have recently settled the Erd\H os distance conjecture in a remarkable paper \cite{GK10}, proving that 
$$ \# \big\{|x-y|: x,y \in P\big\} \gtrsim \frac{n}{\log(n)};$$ from this,  it follows that 

$$ \# T_2(P) \gtrsim \frac{n^2}{\log(n)},$$ 
which, up to logarithmic factors, is the optimal bound. 
However, the main result in \cite{GI10} allowed the second and third listed authors to obtain an upper bound on $u_{2,2}(n)$ for $s$-adaptable sets that is better than the one in (\ref{st}). 

Theorem \ref{kd} will allow us to obtain similar bounds for 
$u_{k,d}(n)$,
the maximal number of times a given $k+1$ point configuration   can arise among a set of $n$ points in ${\Bbb R}^d$, for all $1 \leq k \leq d$, with $d \ge 2$.  We introduce a variant of this quantity:

\begin{definition} \label{almost} Let $P$ be a subset of ${[0,1]}^d$ consisting of $n$ points. For $\delta>0$, define 
$$u^{\delta}_{k,d}(n)=\# \big\{(x^1,x^2, \dots, x^{k+1}) \in P \times P \times \dots \times P: 
t_{ij}-\delta \leq |x^i-x^j| \leq t_{ij}+\delta;1 \leq i<j \leq k+1 \big\},$$ where the dependence on $t=\{t_{ij}\}$ is suppressed. 
\end{definition} 

Observe that obtaining an upper bound for $u^{\delta}_{k,d}(n)$ (with arbitrary $t_{ij}$) immediately
 implies the same upper bound on $u_{k,d}(n)$ defined above. The main result of this section is the following. 

\begin{theorem} \label{discreteresult} Suppose $P \subset {[0,1]}^d$ is $s$-adaptable for $s=d-\frac{d-1}{2k}+a=s_{k,d}+a$ for every sufficiently small $a>0$. Then for every $b>0$, there exists $C_{b}>0$ such that 
\begin{equation} \label{discretekaifest} u^{n^{-\frac{1}{s_{k,d}}-b}}_{k,d}(n) \leq C_b
n^{k+1-\frac{{k+1 \choose 2}}{s_{k,d}}+b}. \end{equation} 

\end{theorem} 

The proof follows from Theorem \ref{kd} in the following way. Let $E$ denote the support of $d\mu^s_P$,  defined as in (\ref{pizdatayamera}) above. We know that if $s>s_{k,d}$, then 
\begin{equation} \label{konchayem} (\mu_P^s \times \mu_P^s \times \dots \times \mu_P^s) \{(x^1,x^2, \dots, x^{k+1}): t_{ij} \leq |x^i-x^j| \leq t_{ij}+\epsilon \} \lesssim \epsilon^{k+1 \choose 2}. \end{equation} 
Taking $\epsilon=n^{-\frac{1}{s}}$, we see that the left hand side of (\ref{konchayem})  is 
$$ \approx n^{-(k+1)} \cdot u_{k,d}^{n^{-\frac{1}{s}}}(n)$$ and we conclude that 
$$ u_{k,d}^{n^{-\frac{1}{s}}}(n) \lesssim n^{k+1-\frac{{k+1 \choose 2}}{s}},$$ which yields the desired result since $s=s_{k,d}+a$.

As we note above, this result is stronger than the previously known $u_{2,2}(n) \lesssim n^{\frac{4}{3}}$. We also see, for instance, that 
$$ u_{2,3}(n) \lesssim n^{\frac{9}{5}}. $$

In the range  $k \ge 2$, $d \ge 3$, to the best of our knowledge no non-trivial estimates for $u_{k,d}(n)$ were previously known.  While our results do not apply to all points sets (recall that they proved  under the additional restriction that $P$ is $s$-adaptable), our conclusion is stronger in that we show not just that any \emph{single} configuration does not repeat very often, but the same holds for the $n^{-\frac{1}{s}}$-\emph{neighborhood} of the configuration.

\vskip.25in 

\subsection{Volumes determined by discrete sets} 

\vskip.125in 

Dumitrescu, Sharir and Toth prove the following result in \cite{DST09}. 
\begin{theorem} \label{DST} Let $P$ be a subset of ${\Bbb R}^d$ of cardinality $n>>1$. Let $area(x,y,z)$ denote the area of the triangle with endpoints $x,y,z$. Then 
\begin{equation} \label{DSTexponent} \# \{(x,y,z) \in P \times P \times P: area(x,y,z)=t \} \leq Cn^{\alpha_d} \end{equation} for any $t \not=0$, with
$$ \alpha_2=\frac{44}{19} \ \text{and} \ \alpha_3=\frac{17}{7}.$$ 
\end{theorem} 

Application of Theorem \ref{volume} and method from the previous subsection yields the following result. 

\vskip.125in 

\begin{theorem} \label{discretevolume} Let $P \subset {[0,1]}^d$, $d \ge 2$, be an $s$-adaptable set for some $s>s_d$, where 
$$ s_d=d-1+\frac{1}{2d} \ \text{when} \ d \ \text{is even}, \ \text{and} \ s_d=d-1+\frac{1}{2(d-1)} \ \text{when} \ d \ \text{is odd}.$$ 

Let $vol_d(x^1,x^2, \dots, x^{d+1})$ denote the volume of the $d$-dimensional simplex with the endpoints 
$x^1, x^2, \dots, x^{d+1}$. Then for $s>s_d$,

\begin{equation} \label{discretevolumeexponents}  \# \{(x^1, \dots, x^{d+1}): t-n^{-\frac{1}{s}} \leq vol_d(x^1, \dots, x^{d+1})=t+n^{-\frac{1}{s}} \} \leq Cn^{d+1-\frac{1}{s}}. \end{equation} 

\end{theorem} 

\vskip.125in 

When $d=2$, the exponent on the right hand side of (\ref{discretevolumeexponents}) equals $3-\frac{4}{5}=\frac{11}{4}$, which is smaller than the exponent $\alpha_2=\frac{44}{19}$ in Theorem \ref{DST} above. Once again, we caution the reader that, although we obtain a better exponent, it is only under the hypothesis of $s$-adaptability. 

In view of Remark \ref{kvolume} and the conversion mechanism of the section, we can prove, more generally, that if $P \subset {[0,1]}^d$, $d \ge 2$, is an $s$-adaptable set for some $s>\frac{d}{2}+\frac{1}{4}$, then 
$$ \# \{(x,y,z): t-n^{-\frac{1}{s}} \leq area(x,y,z) \leq t+n^{-\frac{1}{s}} \} \lesssim n^{3-\frac{1}{s}}.$$ 

As  seen above, in two dimensions this gives a slightly better exponent, in the context of $s$-adaptable sets, than the one obtained in \cite{DST09}, while in three dimensions, our exponent and the one in \cite{DST09} match. In higher dimensions, the results here are the only ones currently known. 

\vskip.125in 

\subsection{Angles determined by discrete sets} 

The results of this subsection are contained in \cite{IMP11}, but are included to indicate the wide applicability of our method. 

The following results were obtained by Pach and Sharir \cite{PS92}, and Apfelbaum and Sharir \cite{AS05}. In \cite{PS92}, it is shown that for a set of $n$ points in $\mathbb{R}^2$, no angle can occur more than $cn^2 \log n$ times. Since there are about $n^3$ triples of points, this implies that there must be at least $c\frac{n}{\log n}$ distinct angles. In \cite{AS05}, it is shown that for a set of $n$ points in $\mathbb{R}^3$, no angle can occur more than $cn^{\frac{7}{3}}$ times, which gives a lower bound of at least $cn^{\frac{2}{3}}$ distinct angles. They also show that for a set of $n$ points in $\mathbb{R}^4$, no angle besides $\frac{\pi}{2}$ can occur more than $cn^{\frac{5}{2}}\beta(n)$ times, where $\beta(n)$ grows extremely slowly with respect to $n$. This means that there must be about $n^{\frac{1}{2}}(\beta(n))^{-1}$ distinct angles.

In dimensions four and higher, no results are currently available. We have the following theorem, which follows from Theorem \ref{mainangle} and the conversion mechanism of this section. 
\begin{theorem} \label{discrete}
Let $P \subset {\mathbb R}^d$, $\#P=n$, $d \ge 2$, be an $s$-adaptable set for $s>\frac{d+1}{2}$. Then 
$$ \# \{(x^1,x^2,x^3) \in P \times P \times P: \theta_0-n^{-\frac{1}{s}} \leq \theta(x^1,x^2,x^3) \leq \theta_0+n^{-\frac{1}{s}} \} \lesssim n^{3-\frac{1}{s}}.$$ 
\end{theorem} 

In dimensions two and three, these exponents are not as good as the results of \cite{AS05,PS92}. However, Theorem \ref{discrete} gives non-trivial exponents in all dimensions. 

\vskip.125in 

We have considered three problems in this section: distribution of simplices, distribution of volumes and distribution of angles. Many other geometric problems can be handled by similar methods. Moreover, combining classical combinatorial techniques with the methods of this section should lead to sharper exponents in many cases. We hope to address these issues in a sequel. 

\vskip.25in


\section{Nontranslation invariant multilinear estimates} 

\vskip.125in

\begin{theorem} \label{generalmultilinear} Define a variable coefficient multilinear operator by 
$$ S_{\mu}(f_1, \dots, f_k)(x)=\int \dots \int f_1(u^1) \dots f_k(u^k) d\mu_x(u^1, \dots, u^k),$$ where $x, u^j \in {\Bbb R}^d$ and $\{\mu_x\}$ are non-negative Borel measures such that the map $x\mapsto \mu_x$ is measurable. Let $\psi$ be a smooth   function supported in the double of the unit ball that is equal to one on the unit ball. Suppose that for some $\gamma_1, \gamma_2>0$ we have
\begin{equation} \label{one}
 \sup_{R>0}
 \int_{\mathbb R^d} \left| \int_{\mathbb R^d} \widehat{\mu}_x(\xi^1, \xi^2, 0 \dots, 0) \psi(x/R) dx \right| d\xi^2 \lesssim {(1+|\xi^1|)}^{-\gamma_1}
  \end{equation} 

\begin{equation} \label{two}
 \sup_{R>0}
  \int_{\mathbb R^d} \left| \int_{\mathbb R^d} \widehat{\mu}_x(\xi^1, \xi^2, 0 \dots, 0) \psi(x/R) dx \right| d\xi^1 \lesssim {(1+|\xi^2|)}^{-\gamma_2}
  \end{equation}

Then, acting on  nonnegative functions, $S_\mu$ is a bounded multilinear form,
$$
 S_{\mu}:L^2_{-\frac{\gamma_1}{2}}({\Bbb R}^d) \times L^2_{-\frac{\gamma_2}{2}}({\Bbb R}^d) \times L^{\infty}({\Bbb R}^d) \times \dots \times L^{\infty}({\Bbb R}^d) \to L^1(\Bbb R^d).
 $$ 
Moreover, the same conclusion follows when \eqref{one} and \eqref{two} are replaced by the more symmetric condition
 \begin{equation} \label{three-sym}
\sup_{R>0}\left[\sup_{\xi^1\in \mathbb R^d} \int_{  \mathbb R^d} {(1+|\xi^1|)}^{\frac{ \gamma_1}{2}}{(1+|\xi^2|)}^{\frac{ \gamma_2}{2}}
 \left|Q_R(\xi^1,\xi^2) \right| d\xi^2 
 \!+\!
 \sup_{\xi^2\in \mathbb R^d} \int_{  \mathbb R^d} {(1+|\xi^1|)}^{\frac{ \gamma_1}{2}}{(1+|\xi^2|)}^{\frac{ \gamma_2}{2}}
 \left|Q_R(\xi^1,\xi^2) \right| d\xi^1 \right]
 <\infty \, , 
  \end{equation}  
  where
  $$
  Q_{R}(\xi^1,\xi^2)= \int_{  \mathbb R^d} \widehat{\mu}_x(\xi^1, \xi^2, 0 \dots, 0) \psi(x/R) dx \, . 
  $$

\end{theorem} 

\begin{remark}
In Theorem \ref{generalmultilinear} (and, as noted previously, in Theorem \ref{translationinvariant}) there is nothing special about the first two coordinates in the assumption  of Fourier decay of the measure. The theorems can be stated in terms of  any two distinct, distinguished coordinates, and correspondingly change the resulting boundedness conclusion, as in the following.
\end{remark}

\subsection{Proof of Theorem \ref{generalmultilinear}} 

\vskip.125in 

We may assume that $f_j$ are nonnegative Schwartz functions. Then since $\psi$ is equal to $1$ on the unit ball, for any $R$ we have
$$ 
{||H_{\mu}(f_1, \dots, f_k)||}_{L^1(B_R)}
\leq \prod_{j=3}^k {||f_j||}_{\infty} \cdot \int \dots \int f_1(u^1) f_2(u^2) d\mu_x(u^1, \dots, u^k)\psi(x/R) dx
$$
$$=\prod_{j=3}^k {||f_j||}_{\infty} \int \int \int \widehat{f}_1(\xi^1) \widehat{f}_2(\xi^2) \widehat{\mu}_x(\xi^1, \xi^2, 0, \dots, 0) \psi(x/R) dx d\xi^1 d\xi^2$$ 
$$ \leq \prod_{j=3}^k {||f_j||}_{\infty} \cdot I \cdot II, $$ where 
$$ I^2=\int \int {|\widehat{f}_1(\xi^1)|}^2 \left| \int \widehat{\mu}_x(\xi^1, \xi^2, 0, \dots, 0) \psi(x/R) dx \right| d\xi^2 d\xi^1$$ and 
$$ {II}^2=\int \int {|\widehat{f}_2(\xi^2)|}^2 \left| \int \widehat{\mu}_x(\xi^1, \xi^2, 0, \dots, 0) \psi(x/R) dx \right| d\xi^1 d\xi^2.$$

By   assumption (\ref{one}), 
$$ I^2 \lesssim \int {|\widehat{f}_1(\xi^1)|}^2 {(1+|\xi^1|)}^{-\gamma_1} d\xi^1,$$ and by   assumption (\ref{two}), 
$$ {II}^2 \lesssim \int {|\widehat{f}_2(\xi^2)|}^2 {(1+|\xi^2|)}^{-\gamma_1} d\xi^2.$$

To obtain the same conclusion of the theorem under   assumption \eqref{three-sym}, we slightly modify the preceding proof. We multiply and divide the integrand in the expression 
$$
\prod_{j=3}^k {||f_j||}_{\infty}
\int \int  \widehat{f}_1(\xi^1) \widehat{f}_2(\xi^2)  Q_R(\xi^1,\xi^2) d\xi^1 d\xi^2
$$
by $(1+|\xi^1|)^{\frac{\gamma_1}{2}}(1+|\xi^2|)^{\frac{\gamma_2}{2}}$. Then we apply the Cauchy-Schwarz inequality with respect to the measure 
$$
(1+|\xi^1|)^{\frac{\gamma_1}{2}}(1+|\xi^2|)^{\frac{\gamma_2}{2}}|Q_R(\xi^1,\xi^2) |\, d\xi^1 d\xi^2
$$
and we use condition \eqref{three-sym} to conclude the proof.
\bigskip

\begin{corollary} \label{hardylittlewoodsobolev} Under the assumptions of Theorem \ref{generalmultilinear},  we have that 
$$ 
S_{\mu}: L^{p_{\gamma_1}}({\Bbb R}^d) \times L^{p_{\gamma_2}}({\Bbb R}^d)  \times L^{\infty}({\Bbb R}^d) \times \dots \times L^{\infty}({\Bbb R}^d) \rightarrow L^1({\Bbb R}^d),
$$ 
for  $\gamma_j <d,\, j=1,2$, and 
\begin{equation} \label{lameexponent} 
p_{\gamma_j}=\frac{2}{1+\frac{\gamma_j}{d}} >1. 
\end{equation}
\end{corollary}

The corollary follows from Theorem \ref{generalmultilinear} by applying the  Hardy-Littlewood-Sobolev embedding of $L^q$ into $L^p_{-s}$, where $s>0$ and $\frac{1}{q}=\frac{1}{p}+\frac{s}{d}$, $1<q<p<\infty$, together with the observation that, if the $L^p$ improving property holds for nonnegative functions, then it holds for all functions.

\begin{corollary}
Suppose that  for $j,\ell\in \{1,\dots , k\}$,  $\Xi_{j,\ell}$ is an ordered $k$-tuple with $\xi^j$ in the $j$th entry, $\xi^\ell$ in the $\ell$th entry and $0$ in the remaining entries. Let   $0<\gamma_j<d$ for $j=1,\dots , k$ and 
$$
  Q_{R}^{j,\ell}(\Xi_{j,\ell})= \int_{  \mathbb R^d} \widehat{\mu} (\Xi_{j,\ell}) \psi(x/R) dx \, . 
  $$
Assume that either for all $j,\ell\in \{1,\dots,k\}$ we have 
\begin{equation} \label{jell}
\sup_{R>0} \int_{  \mathbb R^d} \left|   Q_{R}^{j,\ell}(\Xi_{j,\ell}) \right| d\xi^\ell \lesssim {(1+|\xi^j|)}^{-\gamma_j}\, , 
  \end{equation} 
  or for all $j,\ell\in \{1,\dots,k\}$, $j\neq \ell$, we have 
 \begin{equation} \label{jell2}
\sup_{R>0}\left[\sup_{\xi^j\in \mathbb R^d} \int_{  \mathbb R^d} {(1+|\xi^j|)}^{\frac{ \gamma_j}{2}}{(1+|\xi^\ell|)}^{\frac{ \gamma_\ell}{2}}
 \left|Q_R(\Xi_{j,\ell}) \right| d\xi^ \ell   
 \!+\!
 \sup_{\xi^\ell\in \mathbb R^d} \int_{  \mathbb R^d} {(1+|\xi^j|)}^{\frac{ \gamma_j}{2}}{(1+|\xi^\ell|)}^{\frac{ \gamma_\ell}{2}}
 \left|Q_R(\Xi_{j,\ell}) \right| d\xi^j \right]
 <\infty \, .  
  \end{equation} 
  Then the following Lebesgue-space estimate holds  for $S_\mu$:
$$
S_\mu:\,\, L^{p_1}({\Bbb R}^d)\times \cdots \times L^{p_k}({\Bbb R}^d)\to L^1({\Bbb R}^d),   
$$
for any $1<p_j <2,\, 1\le j\le k$, and $\theta_j$ with $0\le \theta_j\le 1$, $ \sum_{j=1}^k \theta_j=1$,   such that
\begin{equation} \label{jelluu}
\sum_{j=1}^k \frac{1}{p_j} = 1+ \frac{1}{d} \sum_{j=1}^k \theta_j \gamma_j.
 \end{equation} 
\end{corollary}

This  follows from Cor. \ref{hardylittlewoodsobolev} by permuting the placement of the Sobolev spaces in all possible pairs of locations and applying multilinear complex interpolation between Lebesgue spaces.   More precisely, let $\vec E_l$ be the vector having $1$ in the $l$th entry and zero elsewhere. 
 The initial points of the interpolation are  $\frac{1}{ p_{\gamma_j}  }\vec E_j + 
 \frac{1}{p_{\gamma_l}}\vec E_l$, $1\le j\neq l\le k$,  and the intermediate point is
\begin{equation} \label{jellvv}
\Big(\frac 1{p_1}, \dots , \frac{1}{p_k}\Big)=\sum_{\substack{ 1\le j,l\le k\\ j\neq l}} \theta_{j,l} \bigg(\frac{1}{ p_{\gamma_j}  }\vec E_j + 
 \frac{1}{p_{\gamma_l}}\vec E_l\bigg)
  \end{equation} 
for some $0\le \theta_{j,l}=\theta_{l,j}\le 1$ with $\sum_{j\neq l}\theta_{j,l}=1$. 
Relationship \eqref{jelluu} follows from   \eqref{jellvv} using 
\eqref{lameexponent}    setting 
$\theta_j =\frac 12\sum_{s\neq j} \theta_{j,s}$ and noting that $\frac{1}{p_j}=\sum_{s\neq j} \theta_{j,s} \frac{1}{p_{\gamma_j}}$.

\subsection{Sharpness of estimates}

We now show that the bound given by Cor. \ref{hardylittlewoodsobolev} is, in general, sharp. To see this consider the bilinear fractional integration operator 
$$ B_{\gamma}(f,g)(x)=\int \int f(x-u)g(x-v) {({|u|}^2+{|v|}^2)}^{-\frac{2d-\gamma}{2}} dudv.$$ 
(See \cite{Gr92,KSt99,GrK01}  for   more singular operators of   fractional integral  type.)
It is not difficult to check that $B_{\gamma}$ satisfies the assumption of Theorem \ref{generalmultilinear}. Replacing $f(x), g(x)$ by $f(\delta x), g(\delta x)$ and changing variables shows that if 
$$B_{\gamma}: L^p({\Bbb R}^d) \times L^q({\Bbb R}^d) \rightarrow L^r({\Bbb R}^d),$$ then 
$$\frac{1}{p}+\frac{1}{q}-\frac{1}{r} \leq \frac{\gamma}{d}.$$ Plugging in $p=q$ and $r=1$ shows that the conclusion of Cor. \ref{hardylittlewoodsobolev} is, in general sharp. 

One easily checks that the operator $A_2^2$  from (\ref{att}) satisfies the assumption of Theorem \ref{generalmultilinear} with $\gamma=\frac{1}{2}$. This yields $p_{\gamma_j}=\frac{8}{5}$, $j=1,2$ in (\ref{lameexponent}). However, it is not difficult to check that the better bound $L^{\frac{3}{2}}({\Bbb R}^2) \times L^{\frac{3}{2}}({\Bbb R}^2) \to L^1({\Bbb R}^2)$ actually holds. 

This is  in contrast to the situation in the linear case. The direct analog of a $L^p \times L^p \to L^1$ bound in the bilinear case is a $L^p \to L^2$ bound in the linear case. Let $A^d_1f(x)$ be defined as in (\ref{spherical}) above. Since 
$ |\widehat{\sigma}(\xi)| \lesssim {(1+|\xi|)}^{-\frac{d-1}{2}}$ by the method of stationary phase, it follows that 
$$ {||A^d_1f||}_2={||\widehat{A^d_1f}||}_2 \lesssim {\left( \int {|\widehat{f}(\xi)|}^2 {(1+|\xi|)}^{-(d-1)} 
d\xi \right)}^{\frac{1}{2}}$$
$$ \lesssim {||f||}_{L^{\frac{2}{2-\frac{1}{d}}}({\Bbb R}^d)}$$ by the classical Hardy-Littlewood-Sobolev inequality (see, e.g., \cite{St93}). This is precisely the sharp $L^p({\Bbb R}^d) \rightarrow L^2({\Bbb R}^d)$ bound for the spherical averaging operator $A^d_1$, as pointed out following  (\ref{sphericalpq}).

\vskip.25in 


\section{Estimates for multilinear adjoints of translation invariant multilinear generalized Radon transforms}

\vskip.125in 

We will now focus our attention on the translation invariant case where we can write the multilinear generalized Radon transforms as
$$ T_{\mu}(f_1, \dots, f_k)(x)=\int \dots \int f_1(x-u^1) \dots f_k(x-u^k) d\mu(u^1, \dots, u^k)$$
where $\mu$ is a nonnegative Borel measure. Define the $i$-th multilinear adjoint $T_{\mu}^{*_i}$ by
$$ \left\langle T_{\mu}(f_1, \dots f_{i-1}, f_{i}, f_{i+1}, \ldots, f_k), f_{k+1} \right\rangle = \left\langle T_{\mu}^{*_i}(f_1, \dots f_{i-1}, f_{k+1}, f_{i+1}, \ldots, f_k), f_i \right\rangle$$
where $\langle\cdot,\cdot\rangle$ is the $L^2(\mathbb{R}^d)$ inner product.

\begin{theorem}\label{adjoints}\mbox{}
\begin{enumerate}
\item Suppose that 
$$ |\widehat{\mu}(-\xi,\xi, 0, \dots, 0)| \lesssim {(1+|\xi|)}^{-\gamma}$$ for some $\gamma>0$. Then for all $\gamma_1,\gamma_2 > 0$ such that $\gamma = \gamma_1 + \gamma_2$ we obtain the following estimate on nonnegative functions
$$ T_{\mu}^{*_i}:L^2_{-\frac{\gamma_1}{2}}({\Bbb R}^d) \times L^2_{-\frac{\gamma_2}{2}}({\Bbb R}^d) \times L^{\infty}({\Bbb R}^d) \times \dots \times L^{\infty}({\Bbb R}^d)$$
for $i=3,4,\ldots,k$.
\item  Suppose that 
\begin{equation}\label{AdjointsDecay1}
|\widehat{\mu}(0, \xi, 0, \dots, 0)| \lesssim {(1+|\xi|)}^{-\gamma}
\end{equation}
for some $\gamma>0$. Then for all $\gamma_1,\gamma_2 > 0$ such that $\gamma = \gamma_1 + \gamma_2$ we obtain the following estimate on nonnegative functions
$$ T_{\mu}^{*_1}:L^2_{-\frac{\gamma_1}{2}}({\Bbb R}^d) \times L^2_{-\frac{\gamma_2}{2}}({\Bbb R}^d) \times L^{\infty}({\Bbb R}^d) \times \dots \times L^{\infty}({\Bbb R}^d)$$
\item  Suppose that 
\begin{equation}\label{AdjointsDecay2}
|\widehat{\mu}(\xi,0, 0, \dots, 0)| \lesssim {(1+|\xi|)}^{-\gamma}
\end{equation}
for some $\gamma>0$. Then for all $\gamma_1,\gamma_2 > 0$ such that $\gamma = \gamma_1 + \gamma_2$ we obtain the following estimate on nonnegative functions
$$ T_{\mu}^{*_2}:L^2_{-\frac{\gamma_1}{2}}({\Bbb R}^d) \times L^2_{-\frac{\gamma_2}{2}}({\Bbb R}^d) \times L^{\infty}({\Bbb R}^d) \times \dots \times L^{\infty}({\Bbb R}^d)$$
\end{enumerate}
\end{theorem}

\vskip.25in 

To prove Theorem \ref{adjoints}, start by noting
that it is easy to see that $T_{\mu}^{*_i}(f_1, \dots, f_k)(x)$ is equal to
$$ \int \dots \int f_1(x+u^i-u^1) \dots f_{i-1}(x+u^i-u^{i-1}) f_{i}(x+u^i) f_{i+1}(x+u^i-u^{i+1}) \dots f_k(x+u^i-u^k) d\mu(u^1, \dots, u^k) $$
for all $i=1,\ldots,k$. We now proceed to prove all the cases in the theorem. Assume $f_j$ are nonnegative Schwartz functions.

When $i > 2$ we get
$$ 
{||T_{\mu}^{*_i}(f_1, \dots, f_k)||}_{L^1(\mathbb{R}^d)}
\leq \prod_{j=3}^k {||f_j||}_{\infty} \cdot \int \dots \int f_1(x+u^i-u^1) f_2(x+u^i-u^2) dx d\mu(u^1, \dots, u^k) 
$$
$$ = \prod_{j=3}^k {||f_j||}_{\infty} \cdot \int \dots \int f_1(y) f_2(x+u^1-u^2) dy d\mu(u^1, \dots, u^k) $$
and we observe that precisely this quantity came up in the proof of Theorem \ref{translationinvariant}. Since we have the same assumptions as in that theorem then we note that the same proof will work.

When $i=1$ we get
$$ 
{||T_{\mu}^{*_i}(f_1, \dots, f_k)||}_{L^1(\mathbb{R}^d)}
\leq \prod_{j=3}^k {||f_j||}_{\infty} \cdot \int \dots \int f_1(x+u^1) f_2(x+u^1-u^2) dx d\mu(u^1, \dots, u^k) 
$$
$$ = \prod_{j=3}^k {||f_j||}_{\infty} \cdot \int \dots \int f_1(y) f_2(y-u^2) dy d\mu(u^1, \dots, u^k) $$
$$= \prod_{j=3}^k {||f_j||}_{\infty} \cdot \int \dots \int f_1(y) \int \widehat{f}_2(\xi)e^{2\pi i \xi\cdot(y-u^2)}d\xi dy d\mu(u^1, \dots, u^k)$$
$$=\prod_{j=3}^k {||f_j||}_{\infty} \int \widehat{f}_1(-\xi) \widehat{f}_2(\xi) \widehat{\mu}(0, \xi, 0, \dots, 0) d\xi$$
$$ \leq \prod_{j=3}^k {||f_j||}_{\infty} \int \left|\widehat{f}_1(-\xi)\right| \left|\widehat{f}_2(\xi)\right| (1+|\xi|)^{\gamma} d\xi $$
$$ \leq \prod_{j=3}^k {||f_j||}_{\infty} \left( \int {|\widehat{f}_1(\xi)|}^2 {(1+|\xi|)}^{-\gamma_1} d\xi \right) \left( \int {|\widehat{f}_2(\xi)|}^2 {(1+|\xi|)}^{-\gamma_2} d\xi \right) $$
where in the second to last step we used assumption \eqref{AdjointsDecay1}. The case $i=2$ is similar.

\vskip.25in 


\section{Regular value theorem in a fractal setting} 

\vskip.125in 

The regular value theorem in elementary differential geometry says that if $\phi: X \to Y$, where $X$ is a smooth manifold of dimension $n$ and $Y$ is a smooth manifold of dimension $m<n$, then if $\phi$ is a submersion on the set 
$$\{x \in X: \vec{\phi}(x)=y \}, $$ where $y$ is a fixed element of $Y$, then the set 
$$ {\vec{\phi}}^{-1}(y)=\{x \in X: \vec{\phi}(x)=y \}$$ is either empty or is a $(n-m)$-dimensional sub-manifold of $X$. 

In \cite{EIT11}, the authors considered the situation where $Y={\Bbb R}^m$ and $X$ is replaced by $E \times E$, where $E \subset {\Bbb R}^d$ is a set of a given Hausdorff dimension, which, in general, is far from being a smooth manifold. A direct analogue of the regular value theorem would be the  statement that the set 
$$S_{\vec{t}}^{\vec{\phi}}(E_1,E_2)=\{(x,y) \in E_1 \times E_2: \phi_l(x,y)=t_l; 1 \leq l \leq m \}$$ is either empty or has fractal dimension exactly $s_1+s_2-m$, where $s_j$ is the Hausdorff dimension of $E_j$. The examples in \cite{EIT11}, based on arithmetic constructions, show that the lower bound does not in general hold due to the fractal nature of the problem, it is shown in \cite{EIT11} that if $m=1$ and the Monge-Amp\`ere determinant
\begin{equation} det
\begin{pmatrix} 
 0 & \nabla_{x}\phi \\
 -{(\nabla_{y}\phi)}^{T} & \frac{\partial^2 \phi}{dx_i dy_j}
\end{pmatrix}
\neq 0
\end{equation} on the set $\{(x,y): \phi(x,y)=t \}$, then the upper Minkowski dimension of $S_t^{\phi}(E_1,E_2)$ is indeed $\le\, s_1+s_2-1$. 

The multilinear machinery developed in this paper allows us to study the upper Minkowski dimension of the set 
$$ S_{\vec{t},k}^{\Phi}(E_1, \dots, E_{k+1})=\{(x^1, \dots, x^{k+1}) \in E_1 \times E_2 \times \dots \times E_{k+1}: \Phi(x^1, \dots, x^{k+1})=\vec{t} \},$$ where $\Phi$ is defined as in subsection \ref{multidefsubsec} above. 

\vskip.125in 

The techniques in \cite{EIT11} show readily that Theorem \ref{kd} implies the following result. 

\begin{theorem} \label{regvaluetheorem} Let $E_j \subset {[0,1]}^d$ of Hausdorff dimension $s_j$. Under the assumptions of Theorem \nolinebreak\ref{applicationmama},
\begin{equation} \label{dimensioninequality} \overline{dim}_{{\mathcal M}}(S_{\vec{t},k}^{\Phi}(E_1, \dots, E_{k+1})) \leq s_1+s_2+\dots+s_{k+1}-n, \end{equation} 
provided that
$$ s_1+s_2+\dots+s_{k+1}>(k+1) \left(d-\frac{\gamma}{k} \right).$$ 
\end{theorem} 

\begin{remark} The critical exponent provided by Theorem \ref{regvaluetheorem} is known to be sharp in the case $k=1$. In the multilinear case, the issue is sharpness is related to some interesting questions about the distribution of lattice points on varieties of higher co-dimension. This question shall be investigated systematically in a sequel. \end{remark}

\begin{remark} One can  use Theorems \ref{kd} and  \ref{volume} and  \ref{mainangle} 
to provide corollaries of Theorem \ref{regvaluetheorem} in a variety of settings. \end{remark} 

\vskip.25in 



\end{document}